\documentclass[11pt]{article}
\usepackage[margin=1.1in]{geometry}

\usepackage{amsmath, amssymb, amsthm}
\usepackage{mathtools}

\usepackage{graphicx}
\usepackage{hyperref}
\usepackage{algorithm}
\usepackage{algpseudocode}

\usepackage{bbm}
\usepackage{amsmath}
\usepackage{bm}
\usepackage{mathtools}
\usepackage{amssymb}

\usepackage{comment}
\usepackage{enumitem}
\usepackage{caption}
\usepackage{subcaption}
\usepackage{nameref}
\newcommand{\namerefshort}[2]{\hyperref[#1]{#2}}
\usepackage{makecell}

\usepackage{amsopn}

\DeclareMathOperator*{\argmin}{arg\,min}

\def\epsk{{\varepsilon_k}}
\def\Rma{\mathbb{R}^{m}_a}

\newcommand{\headers}[2]{}
\newcommand{\funding}[1]{}
\newcommand{\email}[1]{\texttt{#1}}
\newenvironment{keywords}
  {\par\smallskip\noindent\textbf{Keywords.} }
  {\par\smallskip}

\newenvironment{MSCcodes}
  {\par\smallskip\noindent\textbf{MSC codes.} }
  {\par\smallskip}
\theoremstyle{plain}
\newtheorem{theorem}{Theorem}[section]
\newtheorem{lemma}[theorem]{Lemma}
\newtheorem{proposition}[theorem]{Proposition}
\newtheorem{corollary}{Corollary}
\theoremstyle{definition}
\newtheorem{definition}[theorem]{Definition}

\theoremstyle{remark}
\newtheorem{remark}[theorem]{Remark}

\headers{Time-integrated Optimal Transport}
{Thai P.D. Nguyen, Hong T.M. Chu, Kim-Chuan Toh}

\title{Time-integrated Optimal Transport: A Robust Minimax Framework%
}

\author{
Thai P.D. Nguyen%
\thanks{College of Engineering and Computer Science, VinUniversity
(\email{thai.npd@vinuni.edu.vn}).}
\and
Hong T.M. Chu%
\thanks{(Corresponding author) College of Engineering and Computer Science, VinUniversity
(\email{hong.ctm@vinuni.edu.vn}).}
\and
Kim-Chuan Toh%
\thanks{Department of Mathematics, and Institute of Operations Research and Analytics,
National University of Singapore, Singapore 119076
(\email{mattohkc@nus.edu.sg}).}
}

\date{}

\begin{document}

\maketitle

\begin{abstract}
Comparing time series in a principled manner requires capturing both temporal alignment and distributional similarity of features. Optimal transport (OT) has recently emerged as a powerful tool for this task, but existing OT-based approaches often depend on manually selected balancing parameters and can be computationally intensive. In this work, we introduce the Time-integrated Optimal Transport (TiOT) framework, which integrates temporal and feature components into a unified objective and yields a well-defined metric on the space of probability measures. This metric preserves fundamental properties of the Wasserstein distance, while avoiding the need for parameter tuning. To address the
corresponding computational challenges, we introduce an entropic regularized approximation of TiOT, which can be efficiently solved using a block coordinate descent algorithm. Extensive experiments on both synthetic and real-world time series datasets demonstrate that our approach achieves improved accuracy and stability while maintaining comparable efficiency.
\end{abstract}



\section{Introduction}
\label{sec: intro}
The Optimal Transport (OT) problem, whose history dates back to the seminal works of Monge (1781) \cite{monge1781} and Kantorovich (1942) \cite{kantorovich1942}, has long been a central tool in mathematical analysis and applications. At its core, it concerns the question of how to transport one distribution of mass into another with minimal costs. The ubiquity and versatility of this idea have led itself to links with many branches of mathematics and beyond. Many applications of optimal transport are based on the distance it induces (the Wasserstein distance $\mathcal{W}_p$), which possesses several fundamental properties, such as reproducing the structure of the underlying space in the space of probability measures and metrizing the weak convergence of probability measures. Equipped with these properties, the distance has been used to provide existence, stability, and uniqueness results of solutions of many important PDEs such as Euler flows \cite{Brenier1989eulerflow}, gradient flows \cite{Samtambrogio2015, Otto1998gdflow}. More broadly, this theoretical foundation has provided a rigorous basis for advances in diverse areas, including game theory, economics, statistics \cite{Samtambrogio2015, Ekeland2005}, image processing and shape recognition \cite{gangbo2000, Haker2004}, and even physics \cite{Frisch2002}. 

Recently, optimal transport has also found a natural role in modern machine learning, with applications in areas such as domain adaptation \cite{Courty2017}, computer graphics \cite{Bonneel2016}, and supervised learning \cite{frogner2015}. Underlying these advances are several essential properties of OT. A prominent example is its ability to metrize weak convergence, which has proven crucial in the development of Wasserstein generative adversarial networks \cite{Arjovsky2017}. Another example is OT’s metric properties, which have been exploited to accelerate similarity search algorithms like the Approximating and Elimination Search Algorithm \cite{AESA1992}. In parallel, the introduction of entropic regularization has greatly improved the scalability of OT in large-scale learning problems, giving rise to the Sinkhorn algorithm \cite{CuturiSinkhorn2013}, which combines computational efficiency with a provable linear convergence rate \cite{Franklin1989,  Carlier22}.

Among various data types studied in machine learning, time series play a particularly important role, appearing in healthcare \cite{morid2023time}, finance \cite{Tsay2014}, economics \cite{granger2014}, and industrial processes \cite{Cook2019}. To analyze such data effectively, many machine learning algorithms rely on a well-defined notion of distance to measure the dissimilarity between samples. Dynamic Time Warping (DTW) \cite{dtw1978} has long been regarded as a standard tool for this task. However, the absence of certain desirable properties for example the metric property in some cases yields unexpected inaccuracy \cite{Jain2019_dtwinaccuracy}. Built on a more solid theoretical basis, the optimal transport problem has lately emerged as a viable alternative in comparing the dissimilarity between time series. In order to utilize the Wasserstein distance for time series, one needs to reasonably incorporate the temporal information into the problem. One possible approach is to impose additional time-dependent constraints on the classical OT problem  \cite{Bartl2024AOT, Eckstein2024AOT, Shi2025}. Notably, Adapted Optimal Transport (AOT) \cite{Bartl2024AOT, Eckstein2024AOT}, which considers only bicausal couplings, can preserve several fundamental properties of the classical OT problem. Nonetheless, although this approach is conceptually natural, its empirical effectiveness has not yet been extensively demonstrated. From an orthogonal perspective, one could define the ground metric as a combination of temporal differences and differences in the other dimensions of the time series. This idea has been utilized in \cite{Thorpe2017, TAOT2020}. Compared to AOT, the latter approach is simpler, easier to implement and has shown its promising empirical performance in \cite{TAOT2020}. However, a major disadvantage of this approach is the need for additional parameter tuning to balance the trade-off between the time dimension and other dimensions, which can be costly when dealing with large datasets.

\paragraph{Our contributions}
To address the aforementioned issues of the current OT-based methods for computing distance between time series, our work provides the following contributions: 
\begin{enumerate}[label={-}]
    \item  We propose the Time-integrated Optimal Transport problem (TiOT), a novel minimax framework that automatically balances temporal and feature information, thus eliminating the need for manual parameter tuning. In addition, we prove that TiOT  induces a proper metric, denoted as \(\mathcal{D}_p\), on the space of probability measures, and this metric retains key theoretical properties of the Wasserstein distance $\mathcal{W}_p$, such as metrizing weak convergence.
    
    \item We introduce an entropic regularized counterpart eTiOT as a fast and reliable approximation of the proposed TiOT problem, and, prove that its solutions converge to those of the original TiOT as the regularization parameter tends to zero.
    
    \item We develop a Block Coordinate Descent (BCD) algorithm to solve the eTiOT problem and analyze the convergence of this framework based on the theory of BSUM methods \cite{Razaviyayn2013, Hong2017BCD}. Finally, we demonstrate the empirical effectiveness of the proposed metric \(\mathcal{D}_p\) and the computational efficiency of the BCD algorithm through various numerical experiments, including classification tasks on real-world time series datasets.
\end{enumerate}

Robust minimax optimal transport has been a subject of considerable study; see, e.g., \cite{paty2019,Dhouib2020ROT,Pratik2021, Aquino2020,Huang2024},  and references therein. Nevertheless, these studies do not address the incorporation of temporal information in time series within the optimal transport framework, a setting that, as we later demonstrate, admits strong theoretical guarantees and favorable empirical performance.

\paragraph{Structure of the paper} We begin with notation and preliminaries. Section~\ref{sec: TiOT} introduces the formulations of TiOT and eTiOT and establishes key theoretical results. Section~\ref{sec: BCD} develops a block coordinate descent algorithm for eTiOT and analyzes its convergence. Section~\ref{sec: experiment} presents numerical experiments that validate our analysis and highlight the advantages of our framework in practical scenarios. Finally, Section~\ref{sec: conclusions} summarizes the main findings.

\paragraph{Notation}
We write $P(X)$ for the space of probability measures on $X$. For a measurable map $T : X \to Y$ and a measure $\alpha \in P(X)$, $T_{\#}\alpha$ denotes the pushforward of $\alpha$ through $T$, and $\operatorname{Id}$ is the identity map. The standard inner product on $\mathbb{R}^d$ is denoted as $\langle x,y \rangle = \sum_{i=1}^d x_i y_i$ and its induced Euclidean norm is denoted as $\|x\|_2 = (\sum_{i=1}^d x_i^2)^{1/2}$. The $i$-th canonical basis vector of $\mathbb{R}^d$ is denoted by $e_i = (0,\dots,0,1,0,\dots,0)^\top$, with $1$ in the $i$-th position. The vector of all ones is denoted by 
$\mathbbm{1}_d = (1,1,\dots,1)^\top \in \mathbb{R}^d$.
 Given a weight vector $a = (a_1,\dots,a_d)^\top$ with $a_i>0$ and $\sum_{i=1}^d a_i=1$, the weighted $L^2$ norm is $\|x\|_{L^2(a)} = (\sum_{i=1}^d a_i x_i^2)^{1/2}$. For vectors or matrices $x,y$, $x \circ y$ denotes their element-wise (Hadamard) product, and $x \oslash y$ their element-wise division.
In this paper, which focuses on the analysis of time series data, we work with the underlying space \(\mathbb{R}^{d+1}\), where the last coordinate explicitly represents time. The space of probability 
measures on \((\mathbb{R}^{d+1}, d_p)\) with finite \(p\)-th moment is denoted as 
\begin{equation*}
    \begin{array}{c}
        P_p(\mathbb{R}^{d+1}) := \left \{ \mu \in P_p(\mathbb{R}^{d+1}): \int_{\mathbb{R}^{d+1}} d_p(z_0, z) ^p\mu(dz) < +\infty \right\},
    \end{array} 
\end{equation*}
where $z_0 \in \mathbb{R}^{d+1}$ is arbitrary, $d_p(z, z') = \|z - z'\|_p = \left( \sum_{i = 1}^{d+1} |z_i - z'_i|^p \right)^{1/p}.$

We now recall the definition of the optimal transport problem and the 
associated Wasserstein distance, which serve as the starting point for 
our formulation of the Time-integrated Optimal Transport problem. 
\begin{definition}
    Given $\alpha, \beta  \in P_p(\mathbb{R}^{d+1}),$ and $c: \mathbb{R}^{d+1} \times \mathbb{R}^{d+1} \to [0, +\infty]$, the Optimal Transport problem (OT) is formulated as follows
    \begin{equation}
        \begin{array}{c}
           \inf\left\{\int_{\mathbb{R}^{d+1} \times \mathbb{R}^{d+1}}c(z,z')d\pi(z,z') : \pi \in \Pi(\alpha, \beta)  \right\},
        \end{array}
        \label{def: OT}
    \end{equation}
    where $\Pi(\alpha, \beta) = \{ \pi \in P(\mathbb{R}^{d+1} \times \mathbb{R}^{d+1}): \pi(A \times \mathbb{R}^{d+1}) = \alpha(A), \pi(\mathbb{R}^{d+1} \times B) = \beta(B)  \}$ for any measurable subsets $A,B \subset \mathbb{R}^{d+1}$. 
\end{definition}

 One might regard the optimal value of the OT problem as a way to measure the discrepancy between two measures. In general, it does not satisfy the axioms of distance, but if the probability measures are restricted to the set $P_p(\mathbb{R}^{d+1})$ and the metric $d_p$ is used to construct the cost function, then the formulation (\ref{def: OT}) induces a proper distance, generally known as the Wasserstein distance.
\begin{definition}
    Let $p \in [1,\infty)$ and two probability measures $\alpha, \beta$ in $P_p(\mathbb{R}^{d+1}).$ The Wasserstein distance of order $p$ between $\alpha$ and $\beta$ is defined as follows:
    \begin{equation}
        \begin{array}{c}
            \mathcal{W}_p(\alpha, \beta) = \left( \min_{\pi \in \Pi(\alpha, \beta)} \int_{\mathbb{R}^{d+1} \times \mathbb{R}^{d+1}} d_p(z, z')^p d\pi(z,z')  \right)^{1/p}. 
        \end{array}
        \label{eq: continuousOT}
    \end{equation}
\end{definition}

We conclude this section by recalling the common notion of convergence 
for probability measures, namely weak convergence in \(P_p(\mathbb{R}^{d+1})\).
\begin{definition}[Weak convergence in $P_p(\mathbb{R}^{d+1})$] 
\label{def: weak_convergence}Given $p \in [1, \infty)$, let $(\mu_k)_{k \in \mathbb{N}}$ be a sequence of probability measures in $P_p(\mathbb{R}^{d+1})$ and let $\mu$ be another element of $P_p(\mathbb{R}^{d+1})$. Then $(\mu_k)$ is said to converge weakly in $P_p(\mathbb{R}^{d+1})$ if for any bounded continuous function $\phi $ and $z_0 \in \mathbb{R}^{d+1}$:
\begin{equation*}
    \begin{array}{c}
        \int \phi d\mu_k \to \int \phi d\mu \quad \text{ and } \quad \int d_p(z_0, z)^p d\mu_k(z) \to \int d_p(z_0, z)^pd\mu(z).
    \end{array}
\end{equation*}
     
\end{definition}

\section{Time-integrated Optimal Transport} 
\label{sec: TiOT}
To address the limitations of current OT-based metrics for time series \cite{Bartl2024AOT, Eckstein2024AOT,Shi2025, TAOT2020, Thorpe2017}, we introduce the Time-integrated Optimal Transport problem (TiOT), which explicitly incorporates temporal information into the OT cost matrix. The guiding principle of TiOT is to reformulate the optimal transport problem in a way that achieves the maximum discrimination between two time series. Ideally, this expansion should not be uniform across all pairs of time series: it should be weaker for series that are intrinsically closer to each other and stronger for those that are intrinsically more distant. The formal formulation of this idea is given below, while its numerical advantages will be elaborated in later experiments.
\begin{definition}[Time-integrated Optimal Transport]
    Given $\alpha, \beta  \in P_p(\mathbb{R}^{d+1})$, the Time-integrated Optimal Transport problem between these measures reads
\begin{align}
    \mathcal{D}_p(\alpha,\beta) 
    &= \max\limits_{w \in [0,1]} \left[ \min\limits_{\pi \in \Pi(\alpha, \beta)} 
    \int_{\mathbb{R}^{d+1} \times \mathbb{R}^{d+1}} d_{p,w}((x,t),(y,s))^p \, d\pi((x,t),(y,s)) \right]^{1/p} \label{def: continuousTiOT} \\
    &= \max\limits_{w \in [0,1]} \mathcal{W}_{p,w}(\alpha, \beta), \nonumber
\end{align}
    where $p \geq 1$ and $d_{p,w}((x,t), (y,s)) = \left(w\|x-y\|_p^p + (1-w)|t-s|^p \right)^{1/p}$. In particular, $d_p = 2^{1/p} \times  d_{p,1/2}$ and $\mathcal{W}_p = 2^{1/p} \times \mathcal{W}_{p, 1/2}$.
\end{definition}

\begin{remark}
We give some remarks on the TiOT problem:
\begin{enumerate}[label={-}]
    \item For $\alpha, \beta \in P_p(\mathbb{R}^{d+1})$, the function $d_{p,w}^p$ is integrable with respect to any coupling $\pi$, since $d_p^p$ is integrable with respect to $\pi$ and $d_{p,w}(z,z') \leq d_p(z,z')$ for all $z,z' \in \mathbb{R}^{d+1}$ and $w \in [0,1]$. Thus, the finiteness of $\mathcal{D}_p$ is guaranteed.
    \item In general, $d_{p,w}$ is only a pseudometric, as the positivity axiom may fail when $w=0$ or $w=1$. Nevertheless, we will later show that this is sufficient for $(P_p(\mathbb{R}^{d+1}), \mathcal{D}_p)$ to form a metric space, as a result of the maximum property.
    \item The existence of an optimal solution of (\ref{def: continuousTiOT}) is given in Appendix \ref{appen: proof_existence}.
\end{enumerate}
\end{remark}

In the following sections, we use the notation $\mathcal{W}_{p,w}$ to denote the $p$-Wasserstein distance with respect to the underlying metric $d_{p,w}$, and $\mathcal{W}_p$ to denote the ordinary $p$-Wasserstein distance with the standard metric $d_p$ defined in (\ref{eq: continuousOT}). 

\subsection{Time-integrated Optimal Transport metric space}
In this section, we show that $\mathcal{D}_p$ is a natural extension of the classical Wasserstein distance $\mathcal{W}_p$, in the sense that $\mathcal{D}_p$ preserves several fundamental properties of $\mathcal{W}_p$. We begin by demonstrating that $\mathcal{D}_p$ satisfies the metric property, which is essential for any distance function. We then extend various convergence and topological characterizations from $(P_p(\mathbb{R}^{d+1}), \mathcal{W}_p)$ to $(P_p(\mathbb{R}^{d+1}), \mathcal{D}_p)$.

\begin{theorem}
    The Time-integrated Optimal Transport function $\mathcal{D}_p$ is a distance on $P_p(\mathbb{R}^{d+1}) $.
\end{theorem}
\begin{proof}
We now verify, in sequence, the symmetry, positivity, and triangle inequality for $\mathcal{D}_p$. Let $\alpha, \beta \in P_p(\mathbb{R}^{d+1})$ with supports $X,Y$.  
First, we have $\mathcal{D}_p(\alpha,\beta) = \mathcal{D}_p(\beta,\alpha)$ by the symmetry of $d_{p,w}$ . 

Second, to show $\mathcal{D}_p(\alpha,\alpha)=0$, consider the diagonal coupling 
$\pi_0 := (\operatorname{Id}\times\operatorname{Id})_{\#}\alpha$, that is, $\pi_0(V) = \alpha(\{(x,t): ((x,t),(x,t)) \in V\})$ for any measurable 
$V \subset X\times X$. Then $\pi_0 \in \Pi(\alpha,\alpha)$ and $d_{p,w}((x,t),(x,t))=0$ for all $w \in [0,1]$. Therefore, for any $w \in [0,1]$
\begin{equation*}
\begin{array}{rcl}
\mathcal{W}_{p,w}(\alpha, \alpha) 
&=& \min\limits_{\pi \in \Pi(\alpha, \alpha)}
          \int_{X \times X} d_{p,w}((x,t), (y,s))^p \, d\pi((x,t),(y,s)) \\
&\leq& 
          \int_{X \times X} d_{p,w}((x,t), (y,s))^p \, d\pi_0((x,t),(y,s))  = 0. \\
\end{array}
\end{equation*}
Hence, by definition, \(\mathcal{D}_p(\alpha, \alpha) = \max\limits_{w \in [0,1]} \mathcal{W}_{p,w}(\alpha, \alpha) = 0 \). If $\mathcal{D}_p(\alpha,\beta)=0$, then for any $w \in [0,1]$, $\mathcal{W}_{p,w}(\alpha,\beta)\leq \mathcal{D}_p(\alpha, \beta) = 0$. Thus,  $\mathcal{W}_{p,w}(\alpha,\beta)=0$ for all $w\in[0,1]$.  
Since $\mathcal{W}_{p,w}$ is a valid metric for $w\in(0,1)$, it follows that $\alpha=\beta$.

Third, to prove the triangle inequality, we invoke the Gluing Lemma \cite{villani2003topics}.  
Let $ \xi \in P_p(\mathbb{R}^{d+1})$ be supported on $Z$, and let 
$\pi^*_{xy}, \pi^*_{yz}$ be optimal couplings for $(\alpha,\beta)$ and $(\beta,\xi)$, respectively. The lemma ensures the existence of $\pi \in P(X \times Y \times Z)$ with marginals 
$\pi_{xy}=\pi^*_{xy}$ and $\pi_{yz}=\pi^*_{yz}$. Its marginal $\pi_{xz}$ then belongs to $\Pi(\alpha,\xi)$. Hence,
\begin{equation}\label{eq: TiOT_triangle}
\begin{array}{ll}
\mathcal{D}_p(\alpha, \xi) \!\!\!\!\!
&\leq 
\max\limits_{w \in [0,1]} \left[ 
          \int_{X \times Z} d_{p,w}((x,t), (z,e))^p \, d\pi_{xz} \right]^{1/p}  \\
&= \max\limits_{w \in [0,1]} \left[ 
          \int_{X \times Y \times Z} d_{p,w}((x,t), (z,e))^p \, 
          d\pi \right]^{1/p}  \\
&\leq \max\limits_{w \in [0,1]} \left[ 
          \int_{X \times Y \times Z} 
          \big( d_{p,w}((x,t), (y,s)) + d_{p,w}((y,s), (z,e)) \big)^p 
          d\pi \right]^{1/p},
\end{array}
\end{equation}
where the first inequality comes from the feasibility of $\pi_{xz}$, and the second from the triangle inequality for $d_{p,w}$ (a direct consequence of the $L_p$ norm triangle inequality). Next, we bound the right-hand side (RHS) of (\ref{eq: TiOT_triangle}) 
\begin{equation*}
    \begin{array}{ll}
\text{RHS of } (\ref{eq: TiOT_triangle}) &\leq \max\limits_{w \in [0,1]} \Big\{ 
          \left[ \int_{X \times Y \times Z} d_{p,w}((x,t), (y,s))^p \,
          d\pi \right]^{1/p} 
         \\
&\qquad\qquad 
+ 
          \left[ \int_{X \times Y \times Z} d_{p,w}((y,s), (z,e))^p \,
          d\pi \right]^{1/p} \Big\}  \\
&= \max\limits_{w \in [0,1]} \Big\{ 
          \left[ \int_{X \times Y} d_{p,w}((x,t), (y,s))^p \,
          d\pi_{xy}^* \right]^{1/p} \\
& \qquad\qquad + 
          \left[ \int_{Y \times Z} d_{p,w}((y,s), (z,e))^p \,
          d\pi_{yz}^* \right]^{1/p} \Big\}  \\
&\leq \max\limits_{w \in [0,1]} 
          \left[ \int_{X \times Y} d_{p,w}((x,t), (y,s))^p \,
          d\pi_{xy}^* \right]^{1/p} \\
& \qquad\qquad + 
          \max\limits_{w \in [0,1]} 
          \left[ \int_{Y \times Z} d_{p,w}((y,s), (z,e))^p \,
          d\pi_{yz}^* \right]^{1/p} \\
&= \mathcal{D}_p(\alpha, \beta) + \mathcal{D}_p(\beta, \xi).
\end{array}
\end{equation*}
    Here the first inequality follows from the Minkowski inequality and the last inequality from the property of the max function. This completes the proof. 
\end{proof}

The following proposition establishes that $\mathcal{D}_p$ is equivalent to the classical $\mathcal{W}_p$ distance. This result is crucial, as it enables the natural extension of convergence and topological properties from $(P_p(\mathbb{R}^{d+1}), \mathcal{W}_p)$ to $(P_p(\mathbb{R}^{d+1}), \mathcal{D}_p)$. Ensuring these properties is essential and has been shown to play a significant role in various theoretical and applied contexts \cite{Brenier1989eulerflow, Samtambrogio2015, Otto1998gdflow, Arjovsky2017}.
\begin{proposition}
\label{prop: d_bounds}
    Given $p \in [1, \infty )$, we have
    \begin{equation}
        \left( \frac{1}{2}\right)^{1/p} \mathcal{W}_p \leq \mathcal{D}_p \leq \mathcal{W}_p. \label{eq: d_bounds}
    \end{equation}
\end{proposition}
\begin{proof}
    We begin by proving the first inequality. By the definition of $\mathcal{D}_p$ we have that for any $w \in [0,1]$,  $\mathcal{D}_p \geq \mathcal{W}_{p,w}$.
    Choose $w = \frac{1}{2}$, we obtain for any $\alpha,\beta\in P_p(\mathbb{R}^{d+1})$,
    \begin{equation*}
\begin{array}{rcl}
\mathcal{D}_p (\alpha,\beta)
&\geq& \min\limits_{\pi \in \Pi(\alpha, \beta)} 
        \left( \int_{\mathbb{R}^{d+1} \times \mathbb{R}^{d+1}} 
        d_{p,1/2}((x,t),(y,s))^p \, d\pi \right)^{1/p} \\
&=& \min\limits_{\pi \in \Pi(\alpha, \beta)} 
        \left( \int_{\mathbb{R}^{d+1} \times \mathbb{R}^{d+1}} 
        \tfrac{1}{2}\, d_p((x,t),(y,s))^p \, d\pi \right)^{1/p} =  \left(\tfrac{1}{2}\right)^{1/p} \mathcal{W}_p(\alpha,\beta). 
\end{array}
\end{equation*}
Thus the first inequality is established. The relation $d_{p,w} \leq d_p$ for all $w\in[0,1]$ then yields the second inequality immediately.
\end{proof}

\begin{corollary}
\label{corol: cauchy_tightness}
    Given $p \in [1, \infty)$ and let $(\mu_k)_{k \in \mathbb{N}}$ be a Cauchy sequence in $(P_p(\mathbb{R}^{d+1}), \mathcal{D}_p)$. Then $\{\mu_k\}$ is tight.
\end{corollary}
\begin{proof}
    By the first inequality of (\ref{eq: d_bounds}), if a sequence is a Cauchy sequence in $(P_p(\mathbb{R}^{d+1}), \mathcal{D}_p)$ then it must also be a Cauchy sequence in $(P_p(\mathbb{R}^{d+1}), \mathcal{W}_p)$. Combine with the classical result that Cauchy sequences in $\mathcal{W}_p$ sense are tight \cite[Lemma 6.14]{villani2009oldandnew}, the proof is completed.
\end{proof}

We recall that a natural notion of convergence for measures in $P_p(\mathbb{R}^{d+1})$ is weak convergence (definition \ref{def: weak_convergence}). It is known to be equivalent to convergence induced by the metric $\mathcal{W}_p$, and, by the bounds in (\ref{eq: d_bounds}), we could show that the same property holds for $\mathcal{D}_p$.
\begin{theorem}[convergence in $\mathcal{D}_p$]
\label{thm:weakconvergence_TiOT}
If $(\mu_k)_{k \in \mathbb{N}}$ is a sequence of measures in $P_p(\mathbb{R}^{d+1})$ and $\mu$ is another measure in $P(\mathbb{R}^{d+1})$, then the following two statements are equivalent:
\begin{enumerate}[label=(\roman*), leftmargin=2em]
    \item  $\mu_k \text{ converges weakly in } P_p(\mathbb{R}^{d+1}) \text{ to } \mu$;
    \item $\mathcal{D}_p(\mu_k, \mu) \rightarrow 0$.
\end{enumerate}
\end{theorem}
\begin{proof}
    Assume that $\mu_k  \xrightharpoonup{} \mu$, by 
    \cite[Theorem 6.9]{villani2009oldandnew}, we have $\mathcal{W}_{p}(\mu_k, \mu) \to 0.$
    Since $\mathcal{D}_p \leq \mathcal{W}_p$ by the bound in proposition \ref{prop: d_bounds}, we have $\mathcal{D}_p(\mu_k, \mu) \to 0$ as $\mu_k$ converges weakly to $\mu$.

    Conversely, suppose $\mathcal{D}_p(\mu_k, \mu) \to 0$. Then 
 $0 \leq \left( \frac{1}{2} \right)^{1/p} \mathcal{W}_{p} \leq \mathcal{D}_p$ in proposition \ref{prop: d_bounds} implies that
 $\mathcal{W}_{p}(\mu_k, \mu) \to 0$ . Hence it follows from \cite[Theorem 6.9]{villani2009oldandnew} that $\mu_k$ converges weakly to $\mu$ in $P_p(\mathbb{R}^{d+1})$.
\end{proof}

\begin{corollary}[continuity of $\mathcal{D}_p$] Given $p \in [1, \infty)$,  $\mathcal{D}_p$ is continuous on $P_p(\mathbb{R}^{d+1})$. To be specific, 
if $\mu_k$ (resp. $\nu_k$) converges to $\mu$ (resp. $\nu$) weakly in $P_p(\mathbb{R}^{d+1})$ as $k \to \infty$, then 
\begin{equation*}
    \mathcal{D}_p(\mu_k, \nu_k)  \longrightarrow \mathcal{D}_p(\mu, \nu).
\end{equation*}
\end{corollary}

The corollary on the continuity of $\mathcal{D}_p$ above is a direct consequence of Theorem \ref{thm:weakconvergence_TiOT} and the triangle inequality of $\mathcal{D}_p$. Next we show that $(P_p(\mathbb{R}^{d+1}), \mathcal{D}_p)$ preserves the polish property of the base space $\mathbb{R}^{d+1}$ just like the Wasserstein distance.

\begin{theorem}[topology of $\mathcal{D}_p$ space]
    The space $P_p(\mathbb{R}^{d+1})$ equipped with the TiOT distance $\mathcal{D}_p$ is a complete separable metric space.
\end{theorem}
\begin{proof}
    Since  $(P_p(\mathbb{R}^{d+1}), \mathcal{W}_p)$ is a separable space, there exists a countable dense set $\mathcal{P} \in P_p(\mathbb{R}^{d+1})$ such that for any $\varepsilon > 0$ and $\mu \in P_p(\mathbb{R}^{d+1})$, there exists $\nu \in \mathcal{P}$ such that
        $\mathcal{W}_p(\mu, \nu) \leq \varepsilon,$
    By (\ref{eq: d_bounds}), we have $\mathcal{D}_p(\mu, \nu) \leq \mathcal{W}_p(\mu, \nu) \leq \varepsilon$. Therefore $\mathcal{P}$ is also dense in $(P_p(\mathbb{R}^{d+1}), \mathcal{D}_p)$, thus, $(P_p(\mathbb{R}^{d+1}), \mathcal{D}_p)$ is separable.

    The completeness of $P_p(\mathbb{R}^{d+1})$ with $\mathcal{D}_p$ distance is again a direct consequence of the completeness of $P_p(\mathbb{R}^{d+1})$ with $\mathcal{W}_p$ distance. In fact, let $(\mu_k)_{k\in \mathbb{N}}$ be a Cauchy sequence in $\left(P_p(\mathbb{R}^{d+1}), \mathcal{D}_p\right)$. Since  $\left( \frac{1}{2} \right)^{1/p} \mathcal{W}_p \leq \mathcal{D}_p$ (by the first bound in  (\ref{eq: d_bounds})), so $(\mu_k)_{k \in \mathbb{N}}$ is also a Cauchy sequence in $\left(P_p(\mathbb{R}^{d+1}), \mathcal{W}_p\right)$. Hence, Theorem 6.18 in \cite{villani2009oldandnew} gives us the convergence of $(\mu_k)$ in $\mathcal{W}_p$ sense, which leads to the convergence in $\mathcal{D}_p$ sense by the second bound in (\ref{eq: d_bounds}).
\end{proof}

We now turn to the computational aspects of the TiOT problem. For discrete measures, the TiOT problem reduces to a max–min problem of a bilinear function. By replacing the inner minimization by its dual, this problem can be reformulated as a linear program. The explicit formulation is provided in Appendix \ref{appen: LP_TiOT}.

\begin{definition}[TiOT]\label{def: TiOT}
    The discrete TiOT problem between two discrete measures $\alpha=\sum_{i = 1}^m a_i \delta_{(x_i, t_i)}$ and $ \beta=\sum_{j = 1}^n b_i \delta_{(y_j, s_j)}$ with $x_i, y_j \in \mathbb{R}^d$ is given by
    \begin{align}
    &\max_{w \in [0,1]} \min_{\pi \in \Pi(\alpha, \beta)} \langle C(w) , \pi \rangle,
    \label{eq: TiOT} 
\end{align}
where $c_{i,j}(w) = w\|x_i - y_j\|^p_p + (1-w) |t_i - s_j|^p$ for  $w \in [0,1]$, 
and $\Pi(\alpha,\beta) = \{ \pi \in \mathbb{R}^{m \times n}_+ \,:\, \pi \mathbbm{1}_n = a, \pi^\top \mathbbm{1}_m = b \}.$
\end{definition}

\subsection{Entropic regularized Time-integrated Optimal Transport}
Analogous to the entropic optimal transport, we introduce the entropic regularized Time-integrated Optimal Transport (eTiOT) problem by adding a Kullback--Leibler regularization term to the objective of the TiOT 
problem~\eqref{eq: TiOT}. This term, equivalent to the entropy function $\langle \pi, \log(\pi)\rangle$, is strictly convex and pushes the optimal solution into the interior of the feasible set, thereby simplifying the non-negativity constraint and favoring solutions closer to the independent coupling $\pi = a \circ b$. We show that eTiOT provides a reliable approximation of \namerefshort{def: TiOT}{TiOT} while offering substantial computational advantages, as will be demonstrated in Section~\ref{sec: experiment}.

\begin{definition}[eTiOT]\label{def: eTiOT}
   Let $\alpha=\sum_{i = 1}^m a_i \delta_{(x_i, t_i)}$, $ \beta=\sum_{j = 1}^n b_i \delta_{(y_j, s_j)}$ be two discrete measures with $x_i, y_j \in \mathbb{R}^d$. The formulation of the entropic regularized Time-integrated Optimal Transport problem between them is as follows:
    \begin{align}
    &\max_{w \in [0,1]} \min_{\pi \in \Pi(\alpha, \beta)} \langle C(w) , \pi \rangle + \varepsilon \mathbf{KL}(\pi | a \circ b),     
    \label{eq: eTiOT}
\end{align}
where 
$c_{i,j}(w) = w\|x_i - y_j\|^p + (1-w) |t_i - s_j|^p$ for $w \in [0,1],$
$\Pi(\alpha,\beta) = \{ \pi \in \mathbb{R}^{m \times n}_+ \,:\, \pi \mathbbm{1}_n = a, \pi^\top \mathbbm{1}_m = b \}$, 
 and $\mathbf{KL}(\pi|\gamma) = \sum_{i,j}\pi_{ij}\log(\frac{\pi_{ij}}{\gamma_{ij}}) - \pi_{ij} + \gamma_{ij}.$
\end{definition}

Before presenting an efficient algorithm for solving \namerefshort{def: eTiOT}{eTiOT} in Section \ref{sec: BCD}, we  demonstrate that \namerefshort{def: eTiOT}{eTiOT} is a proper surrogate for the original \namerefshort{def: TiOT}{TiOT} problem: as the regularization parameter \(\epsilon\) decreases to zero, the solution of eTiOT converges to a solution of \namerefshort{def: TiOT}{TiOT}. While this is a standard result in classical optimal transport, (see \cite[Proposition 4.1]{PeyreCuturi2019}), the extension to our framework is not straightforward. The challenges arise from the max-min structure in  TiOT.

\begin{theorem}[convergence with respect to $\varepsilon$]\label{thm:convergence_eps} Consider a positive sequence $\{\varepsilon_k\}$ . We denote $(w_k, \pi^{\epsk}_{w_k})$ as the solution of the \namerefshort{def: eTiOT}{eTiOT} problem with $\varepsilon = \varepsilon_k$ and $\mathcal{S}^*$ as the set of optimal solutions of the \namerefshort{def: TiOT}{TiOT} problem. Then the following statements hold.
\begin{enumerate}[label=(\roman*), leftmargin=2em]
    \item If $\varepsilon_k \xrightarrow{k \to \infty} 0$, 
    , there exists a subsequence of $\{(w_k, \pi^\epsk_{w_k}) \}$ that converges to a point in  $\mathcal{S}^*$.
    \item If $\varepsilon_k \xrightarrow{k \to \infty} \infty$, the whole sequence $\pi^\epsk_{w_k} \xrightarrow{k \to \infty}  ab^\top$.
\end{enumerate}

\end{theorem}
\begin{proof} 
$(i)$ For a fixed value of $w \in [0,1]$, we denote 
\begin{equation*}
    \pi_w = \underset{\pi \in \Pi(\alpha, \beta)}{\argmin}\langle C(w), \pi \rangle, \quad 
    \pi^\varepsilon_w = \underset{\pi \in \Pi(\alpha, \beta)}{\argmin} \, \langle C(w), \pi \rangle + \varepsilon \mathbf{KL}(\pi | \cdot),
\end{equation*}
where, for brevity, we write $\mathbf{KL}(\pi | \cdot)$ to mean $\mathbf{KL}(\pi | a \circ b)$.

Since $[0,1], \Pi(\alpha, \beta)$ are compact, we can extract a subsequence such that (for the sake of simplicity, we keep the same notation) $\{w_{k}, \pi^\epsk_{w_k}\} \to (\widehat{w}, \widehat{\pi})$ and $ \widehat{w} \in [0,1]$, $ \widehat{\pi} \in \Pi(\alpha, \beta)$. Consider an arbitrary optimal solution of the \namerefshort{def: TiOT}{TiOT} problem $(w^*, \pi_{w^*})$.

By the optimality of $(w^*, \pi_{w^*})$ we have 
\begin{align}
    \langle C(w), \pi_{w} \rangle \leq \langle C(w^*), \pi_{w^*} \rangle \leq \langle C(w^*), \pi \rangle \quad   \forall \pi \in \Pi(\alpha, \beta), w \in [0,1],  
    \label{eq: f_minimax_general}
\end{align}
where the first inequality stems from the property  that for any $w \in [0,1]$,
\begin{align*}
    \langle C(w), \pi_{w} \rangle = \underset{\pi \in \Pi(\alpha, \beta)}{\min}\langle C(w), \pi \rangle \leq \max_{w \in [0,1]} \min_{\pi \in \Pi(\alpha, \beta)} \langle C(w) , \pi \rangle = \langle C(w^*), \pi_{w^*} \rangle,
\end{align*}
and the second inequality follows from the property that for all $ \pi \in \Pi(\alpha, \beta),$
\begin{align*}
    \langle C(w^*), \pi_{w^*} \rangle = \max_{w \in [0,1]} \min_{\pi \in \Pi(\alpha, \beta)} \langle C(w) , \pi \rangle = \min_{\pi \in \Pi(\alpha, \beta)} \langle C(w^*) , \pi \rangle \leq \langle C(w^*), \pi \rangle.
\end{align*}
In particular, when setting $w = w_k, \pi = \pi^\epsk_{w^*}$ in (\ref{eq: f_minimax_general}), we get
\begin{align}
    &\langle C(w_k), \pi_{w_k} \rangle \leq \langle C(w^*), \pi_{w^*} \rangle \leq \langle C(w^*), \pi^\epsk_{w^*} \rangle. \label{eq: f_minimax_detail}
\end{align}

Similarly, the optimality of $(w_k, \pi^\epsk_{w_k})$ ensures that for any $w \in [0,1]$ and $\pi \in \Pi(\alpha, \beta)$
\begin{align}
\label{eq:g_minimax_general}
    \langle C(w), \pi^\epsk_w \rangle + \epsk \mathbf{KL}(\pi^\epsk_w | \cdot) \leq \langle C(w_k), \pi_{w_k}^\epsk\rangle + \epsk \mathbf{KL}(\pi_{w_k}^\epsk | \cdot) &\leq 
    \langle C(w_k), \pi \rangle + \epsk \mathbf{KL}(\pi |\cdot).  
\end{align}
Setting $w = w^*, \pi = \pi_{w_k}$ in (\ref{eq:g_minimax_general}),
we get
\begin{align}
    \langle C(w^*), \pi^\epsk_{w^*} \rangle \!+\! \epsk \mathbf{KL}(\pi^\epsk_{w^*} | \cdot) \!\leq \! \langle C(w_k),\pi^\epsk_{w_k} \rangle \!+\! \epsk \mathbf{KL}(\pi^\epsk_{w_k} |\cdot) \! \leq \! \langle C(w_k), \pi_{w_k} \rangle \!+\! \epsk\mathbf{KL}(\pi_{w_k} | \cdot).
\label{eq:g_minimax_detail}
\end{align}
From (\ref{eq:g_minimax_detail}) and (\ref{eq: f_minimax_detail}), we have 
\begin{eqnarray}
 & & \hspace{-4mm}  \epsk \big(\mathbf{KL}(\pi^\epsk_{w^*} |\cdot) \! - \!\mathbf{KL}(\pi^\epsk_{w_k} | \cdot)\big) 
 \nonumber \\
&\leq& \langle C(w_k), \pi^\epsk_{w_k} \rangle - \langle C(w^*), \pi^\epsk_{w^*} \rangle
    \;\leq\;
    \langle C(w_k), \pi^\epsk_{w_k} \rangle 
    \!-\! \langle C(w^*), \pi_{w^*} \rangle 
  \label{eq-squeeze}  \\    
 &\leq&  
 \langle C(w_k), \pi^\epsk_{w_k} \rangle 
    \!-\! \langle C(w_k), \pi_{w_k} \rangle 
  \nonumber  
\;\leq\; 
 \epsk \big(\mathbf{KL}(\pi_{w_k} | \cdot) \!-\! \mathbf{KL}(\pi^\epsk_{w_k} | \cdot)\big).
 \nonumber
\end{eqnarray}
Since $\mathbf{KL}(\pi| \cdot)$ is a continuous function of \(\pi\) and $\Pi(\alpha, \beta)$ is compact,   
$ \mathbf{KL}(\pi^\epsk_{w_k}| \cdot)$ and $\mathbf{KL}(\pi_{w_k}| \cdot)$ are both bounded. Thus as $\epsk \to 0$, by applying the squeeze theorem to \eqref{eq-squeeze}, we get $ 
\underset{k \to \infty}{\lim} \langle C(w_k), \pi^\epsk_{w_k} \rangle = \langle C(w^*), \pi_{w^*} \rangle $, so $\langle C(\widehat{w}), \widehat{\pi} \rangle = \langle C(w^*), \pi_{w^*} \rangle$. Therefore, by  \eqref{eq: f_minimax_general}, we have
\begin{equation}
\label{eq:1st_ineq_minimax}
    \langle C(w), \pi_{w} \rangle \leq \langle C(\widehat{w}), \widehat{\pi} \rangle, \quad \forall w \in [0,1].
\end{equation}
Taking the limit as $\epsk \to 0$ in the second inequality of (\ref{eq:g_minimax_general}), we get 
\begin{equation}
\label{eq:2nd_ineq_minimax}
    \langle C(\widehat{w}), \widehat{\pi} \rangle \leq \langle C(\widehat{w}), \pi \rangle, \quad \forall \pi \in \Pi(\alpha, \beta).
\end{equation}
Combining (\ref{eq:1st_ineq_minimax}) and (\ref{eq:2nd_ineq_minimax}), we have 
\begin{equation}
      \langle C(w), \pi_{w} \rangle \leq  \langle C(\widehat{w}), \widehat{\pi} \rangle \leq \langle C(\widehat{w}), \pi \rangle, \quad \forall w \in [0,1], \pi \in \Pi(\alpha, \beta),
\end{equation}
which yields the optimality of $(\widehat{w}, \widehat{\pi})$
for problem \namerefshort{def: TiOT}{TiOT}. Hence
$(\widehat{w}, \widehat{\pi}) \in \mathcal{S}^*$.

$(ii)$ Similar to the above proof, take the sequence $\{ \varepsilon_k \}$ that tends to $+\infty$ as $k\to \infty$, we have a subsequence $\{w_{k}, \pi^\epsk_{w_k}\} \to (\widehat{w}, \widehat{\pi})$ with $ \widehat{w} \in [0,1]$, $ \widehat{\pi} \in \Pi(\alpha, \beta)$. Using the second inequality of (\ref{eq:g_minimax_general}) with $\pi = a \circ b$, we have 
\begin{equation*}
    \langle C(w_k), \pi_{w_k}^\epsk\rangle + \epsk \mathbf{KL}(\pi_{w_k}^\epsk | a \circ b) \leq \langle C(w_k), a \circ b \rangle + \epsk \times 0, 
\end{equation*}
which implies that 
\begin{equation*}
    \epsk \mathbf{KL}(\pi_{w_k}^\epsk | a \circ b) \leq \langle C(w_k), a \circ b \rangle \leq 
    \sum_{i,j} |C(w)_{ij}| a_i b_j
    \leq \|C(w)\|_\infty
    \sum_{i,j}  a_i b_j 
    \leq  \|C\|_\infty,
\end{equation*}
where  $\|C\|_\infty = \max_{w \in [0,1]} \{ \|C(w)\|_\infty
\}$, and $\|C(w)\|_\infty$ is the maximum absolute value among 
all the entries of $C(w)$.
Dividing both sides by $\epsk$ and letting this value tends to $\infty$, we get
\begin{equation*}
   \mathbf{KL}(\pi_{w_k}^\epsk | a \circ b) \to 0.
\end{equation*}
Thus, by the Pinsker's inequality 
\cite{Pinsker1964}, $\pi^\epsk_{w_k} \to a \circ b$. Then, the compactness of $\Pi(\alpha, \beta)$ ensures the convergence to $a \circ b$ of the whole sequence.
\end{proof}

\section{Block coordiate descent algorithm for solving \namerefshort{def: eTiOT}{eTiOT}}\label{sec: BCD}

One of the key developments that has substantially advanced the field of optimal transport is the introduction of its entropic regularized formulation \cite{CuturiSinkhorn2013, Benamou2015}. The Sinkhorn algorithm provides an efficient method for solving this problem, thereby enabling the application of optimal transport to large-scale machine learning tasks. 

Indeed, the Sinkhorn algorithm can be interpreted as a block coordinate descent method applied to the dual formulation of the entropic OT problem. Motivated by this perspective, in this section we introduce a block coordinate descent algorithm for the \namerefshort{def: eTiOT}{eTiOT} problem and establish its convergence under the natural normalization proposed in \cite{Carlier22}.

The Lagrangian of the inner minimization problem is given by
\begin{equation*}
   \mathbbm{L}(\pi, u, v) = \langle \pi, C(w) \rangle + \varepsilon \mathbf{KL}(\pi| a \circ b) + u^T(a - \pi \mathbbm{1}_n ) + v^T(b - \pi^T\mathbbm{1}_m).
\end{equation*}
Setting its gradient  \(\frac{\partial \mathbb{L}}{\partial \pi_{ij}}  
= C_{ij}(w) + \varepsilon \log\Big(\frac{\pi_{ij}}{a_i b_j}\Big) - u_i - v_j = 0\), we get
\begin{equation*}
\begin{array}{c}
    \pi_{ij} = a_i b_j \exp\Big(\frac{u_i + v_j - c_{ij}(w)}{\varepsilon}\Big) \quad \forall\; i\in[m],\; j\in[n].
\end{array}
\end{equation*}
Combining with the normalization in \cite{Carlier22}, the   \namerefshort{def: eTiOT}{eTiOT} problem is now given by 
\begin{equation}
    \begin{array}{c}
        \min_{\substack{w \in [ 0 , 1], \\ u\in \Rma, v  \in \mathbb{R}^n}} F(u,v,w) =  -u^\top a - v^\top b + \varepsilon \sum_{i, j = 1}^n \exp \left( \frac{u_i + v_j - c_{ij}(w)}{\varepsilon} \right)a_i b_j - \varepsilon,
    \end{array} \label{eq: dual_eTiOT}
\end{equation}
where $\Rma = \{u \in \mathbb{R}^m : a^\top u = 0\}$.

We apply the block coordinate descent method to solve (\ref{eq: dual_eTiOT}) with exact minimization for the two blocks $u,v$.
For the block $w$, we use one step of the projected gradient descent method. This creates a hybrid combination of the classic BCD algorithm and the block coordinate gradient projection (BCGP) in \cite{Beck13_BCGD}. The updating scheme of the hybrid block coordinate descent algorithm (HBCD) is given as follows:
\begin{equation}
\begin{array}{lllll}
u^{k+1} & = &  \argmin_{u \in \Rma}\varphi^k_u (u) &=& - \varepsilon \log\Big(\exp(\frac{-C(w^k)}{\varepsilon}) (\exp(v^k/\varepsilon) \circ b)\Big) + \lambda^k \mathbbm{1}, \\[1mm]
v^{k+1} & =&  \argmin_{v \in \mathbb{R}^{n}}\varphi^k_v (v) &=& - \varepsilon \log\Big(\exp(\frac{-C(w^k)^\top}{\varepsilon}) (\exp(u^{k+1}/\varepsilon) \circ a)\Big), \\[1mm]
w^{k+1} & =&  \argmin_{w \in [0,1]}\varphi^k_w (w) &=& \operatorname{Proj}_{[0,1]}\Big( w^k - \eta \nabla_w F(u^{k+1}, v^{k+1}, w^k)\Big),
\end{array}
\label{eq: HBCD}
\end{equation}
where the upper-bound functions associated with the blocks $u,v,w$ are
\begin{equation}
\label{eq: upperbounds}
\begin{array}{lll}
   \varphi_u^k(u) &=&  F(u,v^k,w^k), \\
    \varphi_v^k(v) &=& F(u^{k+1},v,w^k), \\
   \varphi_w^k(w) &=&  F(u^{k+1}, v^{k+1}, w^k) 
    + \nabla_w F(u^{k+1}, v^{k+1}, w^k)\,(w - w^k) 
    + \tfrac{1}{2\eta}(w - w^k)^2,
\end{array}
\end{equation}
with $\lambda^k = \varepsilon a^\top \log(\exp(\frac{-C(w^k)}{\varepsilon})(\exp(v^k/\varepsilon) \circ b))$ as the normalizing constant that enforces the extra constraint $a^\top u = 0$ and $\eta$ is the stepsize.

If the stepsize satisfies $\eta \le 1/L_w$, with $L_w$ denoting the  Lipschitz constant corresponding to the block $w$ of $\nabla_w F$ \cite[eq.~(2.3)]{Beck13_BCGD}, then $\varphi_w^k(w)$ is a valid upper-bound function. In this case, the proposed HBCD algorithm belongs to the class of block successive upper-bound minimization (BSUM) methods \cite{Razaviyayn2013, Hong2017BCD}, which address minimization problems by iteratively minimizing block-wise surrogate upper bounds of the objective. Based on \cite{Carlier22, Beck13_BCGD, Razaviyayn2013, Hong2017BCD}, we derive convergence guarantees for HBCD.
 
Lemma \ref{lem: boundedness} and lemma \ref{lem: sufficient_decrease} 
below follow directly from the results in \cite{Carlier22} and \cite{Beck13_BCGD}. For completeness, we include their full proofs in the Appendix \ref{appen: proof_bounded}, \ref{appen: proof_sufficient_decrease}.

\begin{lemma}[boundedness of iterates]
\label{lem: boundedness}
    For every $k \geq 1$, the iterates generated by \eqref{eq: HBCD} satisfy the bounds
    \begin{align*}
        \|u^k\|_\infty \leq 2 \|C\|_\infty, \quad
        \|v^k\|_\infty \leq 3 \|C\|_\infty,
    \end{align*}
    where $\|C\|_\infty = \max_{w \in [0,1]} \big\{ \|C(w)\|_\infty \big\} = \max\{ \|C(0)\|_\infty, \|C(1)\|_\infty \}$.
\end{lemma}

\begin{lemma}[block Lipschitz continuity of $\nabla_w F$] Let $\xi^k \coloneqq (u^k, v^k, w^k)$ be the sequence generated by \eqref{eq: HBCD}. For any $k \geq 1$ and $(w, w') \in [0,1]^2$, we have
\label{lem: block_lipschitz}
\begin{equation}
\label{eq: block_lipschitz}
    |\nabla_w F(u^k, v^k, w) - \nabla_w F(u^k, v^k, w') | \leq (\|\Tilde{C}\|^2_\infty / \varepsilon )\exp \left(\frac{6\|C\|_\infty}{\varepsilon}\right) |w - w'|,
\end{equation}
where $\|\Tilde{C}\|_\infty = \max\{ \large| ||x_i - y_j ||^p_p - |t_i - s_j|^p \large| : i \in [m], j \in [n]\}$.
\end{lemma}
\begin{proof}
    See Appendix \ref{appen: proof_block_lipschitz}.
\end{proof}

\begin{lemma}[sufficient descent property]
\label{lem: sufficient_decrease}
    Let $\zeta^k \coloneqq (u^k, v^k, w^k)$ and $\zeta^{k+1}\coloneqq (u^{k+1}, v^{k+1}, w^{k+1}) $ be generated by \eqref{eq: HBCD} with $\eta = (\varepsilon / \|\Tilde{C}\|^2_\infty) \exp({\frac{-6\|C\|_\infty}{\varepsilon}})$, then the following inequality holds:
    \begin{equation*}
       \begin{array}{rl} 
        F(\zeta^k) - F(\zeta^{k+1}) &\geq \kappa \left( \|u^k - u^{k+1}\|^2_{L^2(a)} \! + \! \|v^k - v^{k+1}\|^2_{L^2(b)} \right) 
         + \tau|w^k - w^{k+1}|^2,
    \end{array}
    \end{equation*}
where  $\kappa = \exp{ (\frac{-6\|C\|_\infty}{\varepsilon})} / 2\varepsilon$, $\tau =  \|\Tilde{C}\|_\infty^2 \exp{ (\frac{ 6\|C\|_\infty}{\varepsilon})} / 2\varepsilon$.
In the above, $ \|u^k - u^{k+1}\|^2_{L^2(a)} = \sum_{i=1}^m 
(u^{k+1}_i-u^k_i)^2 a_i$, and $\|v^k - v^{k+1}\|^2_{L^2(b)}$
is similarly defined. 
\end{lemma}

Relying on the preceding lemmas, we invoke \cite{Razaviyayn2013} to establish the asymptotic convergence of the iterations \eqref{eq: HBCD} and employ the techniques of \cite{Hong2017BCD}  to obtain the sublinear convergence rate.
\begin{theorem}  \label{thm: convergence_xi}
        Let $\Tilde{\mathcal{S}}$ be the optimal solution set of (\ref{eq: dual_eTiOT}) and $\{\xi^k\}$ be the sequence generated by the \eqref{eq: HBCD} with $\eta = (\varepsilon / \|\Tilde{C}\|^2_\infty) \exp({\frac{-6\|C\|_\infty}{\varepsilon}})$. Then $\xi^k$ converges to $\Tilde{\mathcal{S}}$
       in the sense that
        \[
\lim_{k\to\infty}\inf\nolimits_{\xi\in\tilde{\mathcal S}}
\|\xi^k-\xi\|_2 = 0.
\]
\end{theorem}
\begin{proof}
To prove the convergence of HBCD, it is sufficient to show that HBCD iterates \eqref{eq: HBCD} and the objective function $F(.)$ in \eqref{eq: dual_eTiOT} satisfy the assumption in \cite[theorem 2-b]{Razaviyayn2013}.

First, Lemma~\ref{lem: boundedness} implies that the iterates generated by HBCD lie in a compact set.

Second, the upper-bound functions $\varphi_u^k$ and $\varphi_v^k$ in \eqref{eq: upperbounds} satisfy \cite[Assumption~2]{Razaviyayn2013}. Moreover, by lemma \ref{lem: block_lipschitz} and the block descent lemma \cite[Lemma~ 3.2]{Beck13_BCGD}, the function $\varphi_w^k$ also satisfies this assumption.

Third, the strict convexity of the upper-bound functions guarantees that their corresponding minimization subproblems admit unique solutions. Moreover, the differentiability of $F(\xi)$ ensures the regularity condition. Hence, by \cite[Theorem~2-b]{Razaviyayn2013}, the sequence ${\{\xi^k\}}$ generated by HBCD \eqref{eq: HBCD} converges to the set of stationary points, which in our case is the set of optimal solution due to the convexity of $F(\cdot)$. 
\end{proof}

\begin{lemma}
\label{lem: solution_boundedness}
    Let $\mathcal{\Tilde{S}}$ be the optimal solutions set of \eqref{eq: dual_eTiOT}. There exists $\xi^* = (u^*, v^*, w^*) \in \mathcal{\Tilde{S}}$ such that $\|u^*\|_\infty \leq 2\|C\|_\infty$ and $\|v^*\|_\infty \leq 3\|C\|_\infty$.
\end{lemma}
\begin{proof}
 See Appendix \ref{appen: proof_solution_boundedness}.
\end{proof}
\begin{theorem}
\label{thm: sublinear}
    Let $\zeta^k = (u^k, v^k, w^k)$ be the  sequence generated by the iterations \eqref{eq: HBCD} with $\eta = (\varepsilon / \|\Tilde{C}\|^2_\infty) \exp({\frac{-6\|C\|_\infty}{\varepsilon}})$. 
    For any $k \geq 1$, we have
    \begin{equation}
        F(\zeta^k) - F^* \leq \frac{\rho_1 \rho_2}{k},
    \end{equation}
    where $F^* \!$ denotes the optimal value of \eqref{eq: dual_eTiOT}, $\rho_1 \!\!=\!\! (192m + 216n + 24)\frac{\|C\|^2_\infty}{\varepsilon}\exp(\frac{18\|C\|_\infty}{\varepsilon}) $, $\rho_2 = \max \{4/\rho_1 -2 , F(\zeta^1) - F^*, 2 \}$.
\end{theorem}
\begin{proof}
    See Appendix \ref{appen: proof_sublinear}.
\end{proof}

We conclude this section by reformulating HBCD into an efficient and easily implementable algorithm (Algorithm \ref{alg:HBCD}), with several remarks: 
\begin{itemize}[label={-}]
\item To reduce computational cost, define \(g = a \circ \exp(u/\varepsilon)\), \(h = b \circ \exp(v/\varepsilon)\), \(K(w) = \exp(-C(w)/\varepsilon)\). Then the  
HBCD iteration can be simplified to a Sinkhorn-like algorithm. 
\item Algorithm~\ref{alg:HBCD} provides a competitive framework for solving the \namerefshort{def: eTiOT}{eTiOT} problem. By skipping the normalization step (\ref{eq: HBCD}), and updating $w$ and the stopping criterion only once every \texttt{freq} iterations, it avoids redundant operations and reduces computation without compromising the convergence.
\end{itemize}

\begin{algorithm}[h]
\caption{HBCD algorithm for solving  \namerefshort{def: eTiOT}{eTiOT} problem \eqref{eq: eTiOT}.}\label{alg:HBCD}
\begin{algorithmic}[1]
\State \textbf{Input} two discrete distributions  $\alpha=\sum_{i = 1}^m a_i \delta_{(x_i, t_i)}$ and $ \beta=\sum_{j = 1}^n b_i \delta_{(y_j, s_j)}$, entropic regularization parameter $\varepsilon > 0$, $\eta$,  freq.

\State \textbf{Initialize}
$\Gamma \in \mathbb{R}^{m \times n}$ defined by $\Gamma_{ij} = \|x_i-y_j\|_2^2,$
\State \hspace{4.4em} $\Phi \in \mathbb{R}^{m \times n}$ defined by $\Phi_{ij} = |t_i-s_j|^2,$ 
\State \hspace{4.4em} $w = 0.5, \; C = w \Gamma + (1-w)  \Phi, \; K = \exp\!\big(-C/\varepsilon\big)$,\; 
$h = \tfrac{1}{m}\mathbbm{1}_m,\;$

\While{termination criteria not met}
    \State $g \gets a \oslash (K h)$
    \State $h \gets b \oslash (K^\top g)$
    \If{$\text{mod}( t, \text{freq})  =  0$}
        \State $w \gets \max \big\{ \min \big\{ w - \eta \left[ g^\top ((\Phi - \Gamma) \circ K) h \right], 1 \big\}, 0\big\}$ 
        \State $C \gets w \, \Gamma + (1-w)\, \Phi$
        \State $K \gets \exp(-C/\varepsilon)$
    \EndIf    
\EndWhile
\State $\pi \gets \operatorname{Diag}(g) \cdot K \cdot \operatorname{Diag}(h)$,\, $\mathcal{D}_p^{\varepsilon}(\alpha,\beta)  \gets \langle C, \pi \rangle$
\State \Return $\pi, w, \mathcal{D}_p^{\varepsilon}(\alpha,\beta)$ 
\end{algorithmic}
\end{algorithm}

\section{Experiments}
\label{sec: experiment}In this section, we evaluate the empirical effectiveness of the proposed Time-integrated Optimal Transport (\namerefshort{def: TiOT}{TiOT}) problem through extensive numerical experiments. Specifically, we analyze the stability and reliability of the optimal solutions of the \namerefshort{def: TiOT}{TiOT} problem in Section~\ref{subsec: alignment_experiment}; examine the properties of the induced distance in Section~\ref{subsec: lag_experiment}, validate the theoretical convergence and demonstrate the algorithmic benefits of the entropic variant (\namerefshort{def: eTiOT}{eTiOT}) in  Section~\ref{subsec: runtime_experiment}; and finally, we apply TiOT to time series classification problem on several standard datasets in Section~\ref{subsec: knn_experiment}. All experiments are implemented in Python and executed on a machine equipped with a 12th Gen Intel(R) Core(TM) i7-1260P 2.10 GHz processor. The source code is available at \url{https://github.com/Thai-npd/TiOT-code}

For all experiments in this section, we consider distributions 
$\alpha=\sum_{i=1}^m a_i \delta_{(x_i,t_i)}$ and $\beta=\sum_{j=1}^n b_j \delta_{(y_j,s_j)}$ 
with uniform weights, $a=\tfrac{1}{m}\mathbbm{1}_m$ and $b=
\tfrac{1}{n}\mathbbm{1}_n$, 
where $\{x_i\}$ and $\{t_i\}$ are standardized using Z-score normalization; Similarly for $\{y_j\}$ and $\{s_j\}$.

In Sections \ref{subsec: runtime_experiment}, \ref{subsec: knn_experiment}, the eTiOT problem is solved by Algorithm \ref{alg:HBCD} with termination criterion
$\|g \circ (Kh) - a\|_1 < 0.005$, adopted from \cite{TAOT2020} for fair comparisons. To improve efficiency, multiple subiterations are performed for the $w$-block in each iteration: Projected Gradient Descent is repeated until successive changes in $w$ fall below $10^{-7}$ for small-scale problems $(n,m<1000)$ or $10^{-2}$ for large-scale ones, with at most 50 subiterations. The stepsize is set as $\eta = \sigma/20$ if $\sigma \geq 10$ and $\eta = \sigma/10$ otherwise, where $\sigma = \tfrac{1}{\varepsilon}\big(g^\top \big((\Phi - \Gamma)^2 \circ K\big)h\big)$ approximates the local curvature constant in $w$ of the objective minimized by Algorithm~\ref{alg:HBCD}.

\subsection{Robustness of TiOT}
\label{subsec: alignment_experiment}
In this section, we illustrate the robustness of  the transportation plan \(\pi\) generated by our proposed model \namerefshort{def: TiOT}{TiOT} compared to the transportation plan induced by the usual Wasserstein distance with the base metric $d_{2,w}$, that is, $\mathcal{W}_{2,w}$.

\begin{figure}[htb]
    \centering
        \begin{subfigure}[b]{0.325\textwidth}
        \includegraphics[width=\textwidth]{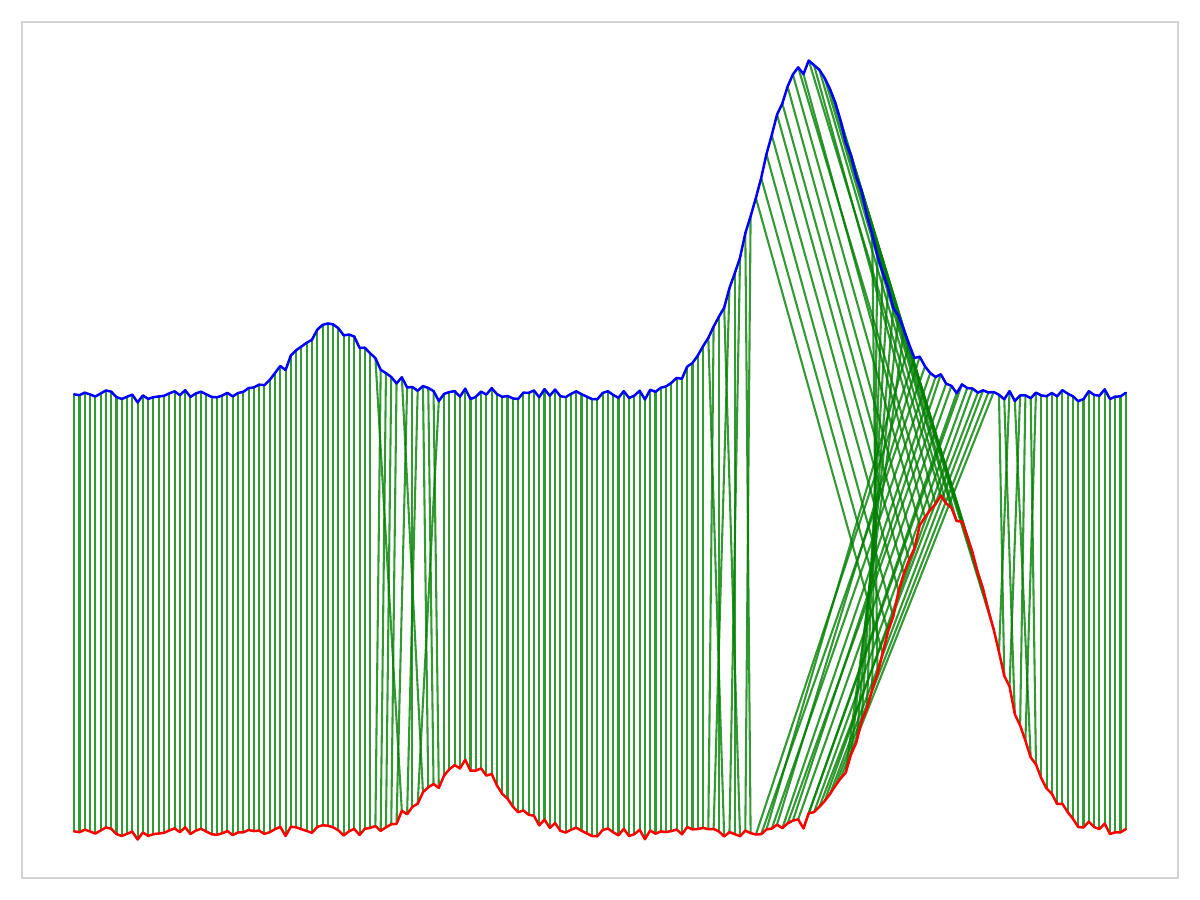}
        \caption{$w = 0.1$}
    \end{subfigure}
    \begin{subfigure}[b]{0.325\textwidth}
        \includegraphics[width=\textwidth]{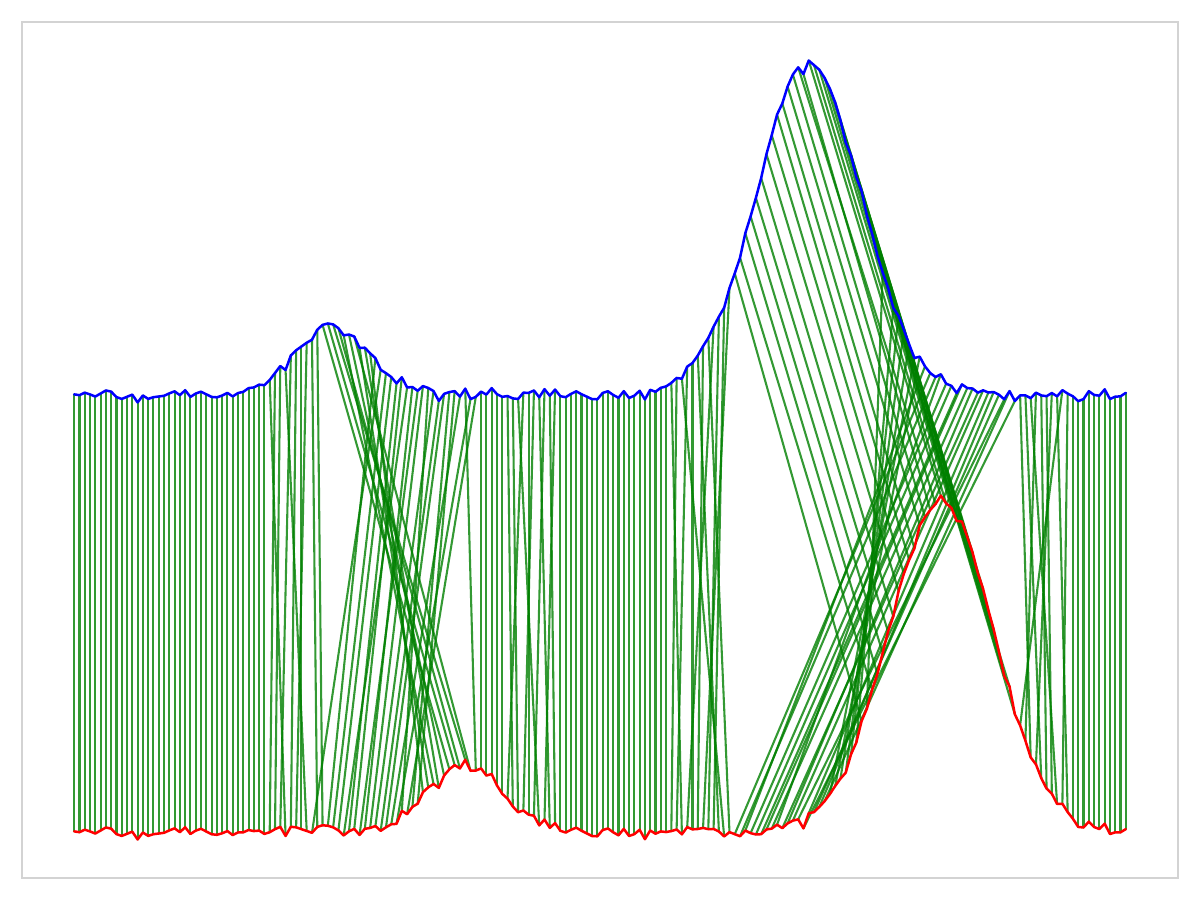}
        \caption{$w = w^* = 0.2867$}
    \end{subfigure}
        \begin{subfigure}[b]{0.325\textwidth}
        \includegraphics[width=\textwidth]{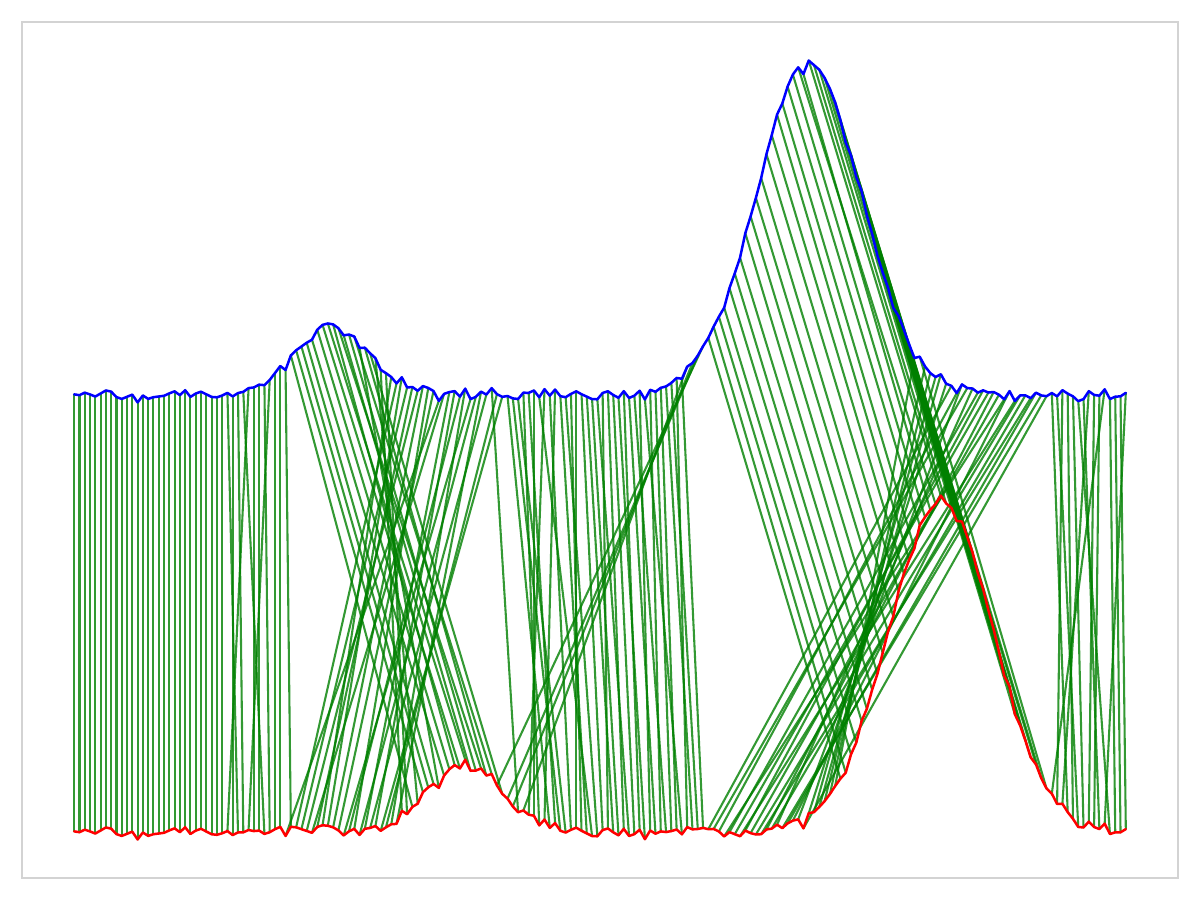}
        \caption{$w = 0.8$}
    \end{subfigure}
    \vspace{-10pt}
    \caption{Alignment (transportation) map \(\pi\) between two time series. Left and right: \(\pi\) induced by $\mathcal{W}_{2,w}$ with selected \(w=0.1\) and \(w=0.8\), middle: \(\pi^*\) induced by our proposed TiOT model.}
        \label{fig:alignment}
\end{figure}

\paragraph{Dataset and experiment design} In this experiment, we generate two mixture of Gaussians time series $\alpha=\sum_{i = 1}^m a_i \delta_{(x_i, t_i)}$ and $ \beta=\sum_{j = 1}^n b_i \delta_{(y_j, s_j)}$ by letting \(m=n=200\),  \(a_i = b_i = 1/200 \), and  \(t_i = s_i = i \)   for any \(i=1,\dots,200\). In addition, we generate \(x_i\)  and \(y_j\) by 
\begin{equation}
\begin{array}{rl}   
    x_i &= 0.2 \exp\left(-\frac{(t_i - 50)^2}{ 2\times 7^2}\right) +  \exp\left(-\frac{(t_i - 140)^2}{2 \times 10^2}\right) + \mathcal{N}(0, 0.01^2), \\
    y_i &=0.2 \exp\left(-\frac{(t_i - 75)^2}{2 \times 7^2}\right) + \exp\left(-\frac{(t_i - 165)^2}{2 \times 10^2}\right) + \mathcal{N}(0, 0.01^2), 
\end{array}
\label{eq: gauss_series}
\end{equation}
for any \(i=1,\dots,200\) and $\mathcal{N}(0, 0.01^2)$ is Gaussian noise. In Figure~\ref{fig:alignment}, we visualize these two mixtures of Gaussians in blue and red, respectively. Finally, we compute \(\mathcal{W}_{2,w} \) with \(w=0.1\) and \(0.8\) to obtain the transportation plan \(\pi\), and solve our proposed TiOT problem to get both optimal parameter \(w^*\) and transportation plan \(\pi^*\). In Figure~\ref{fig:alignment}, each green line represents a non-zero entry \(\pi_{i,j}\).

\paragraph{Analysis} 
We first note that all three alignments are indeed one-to-one alignments, a finding consistent with the fact that one of optimal plans between two discrete distributions of the same size and uniform weight must be a permutation \cite{PeyreCuturi2019}. However, the resulting alignment changes significantly when \(w\) changes. For instance, the temporal constraint imposed by $\mathcal{W}_{2, w=0.1}$ proves overly rigid, preventing the correct alignment of the two peaks. Conversely, $\mathcal{W}_{2,w= 0.8}$ is too permissive, leading to multiple mismatches. In contrast, TiOT adaptively selects $w$ to maximize $\mathcal{W}_{2,w}$, thereby discouraging pairings that might otherwise appear optimal under a poorly chosen weight. This allows TiOT to avoid overfitting and yield a more robust and balanced alignment.

\subsection{Time Series Analysis with $\mathcal{D}_2$ and $\mathcal{W}_{2,w}$ metrics}
\label{subsec: lag_experiment}
In this section, we compare the proposed metric $\mathcal{D}_2$ with $\mathcal{W}_{2,w}$ on a time series lag analysis task. While 
the performance of $\mathcal{W}_{2,w}$ is highly sensitive to its parameter \(w\), $\mathcal{D}_2$ avoids this limitation through its auto-selection mechanism, yielding a more robust and reliable metric.

\paragraph{Dataset and experiment design} The dataset contains daily temperature of Delhi, India, from January 1, 2013, to January 1, 2017 \cite{dailyclimate}, denoted by \(x = (x_1, \dots,x_{1462}) \in \mathbb{R}^{1462}\). A one-year time series starting on day \(\ell\) is defined as
\begin{equation*}
    x^{(\ell)} \coloneqq (x_{\ell},x_{\ell+1},\dots,x_{\ell+364}) \in \mathbb{R}^{365},
\end{equation*}
for \(\ell=1,2,\dots,730\). We then measure the dissimilarity between the initial series \(x^{(1)}\) and its lagged versions \(x^{(\ell)}\) using $\mathcal{D}_2$ \eqref{def: TiOT} and $\mathcal{W}_{2,w}$ with fixed parameters \(w=0.2, 0.5, 0.8\).

\begin{figure}[htb]
    \centering
        \begin{subfigure}[b]{0.495\textwidth}
        \includegraphics[width=\textwidth]{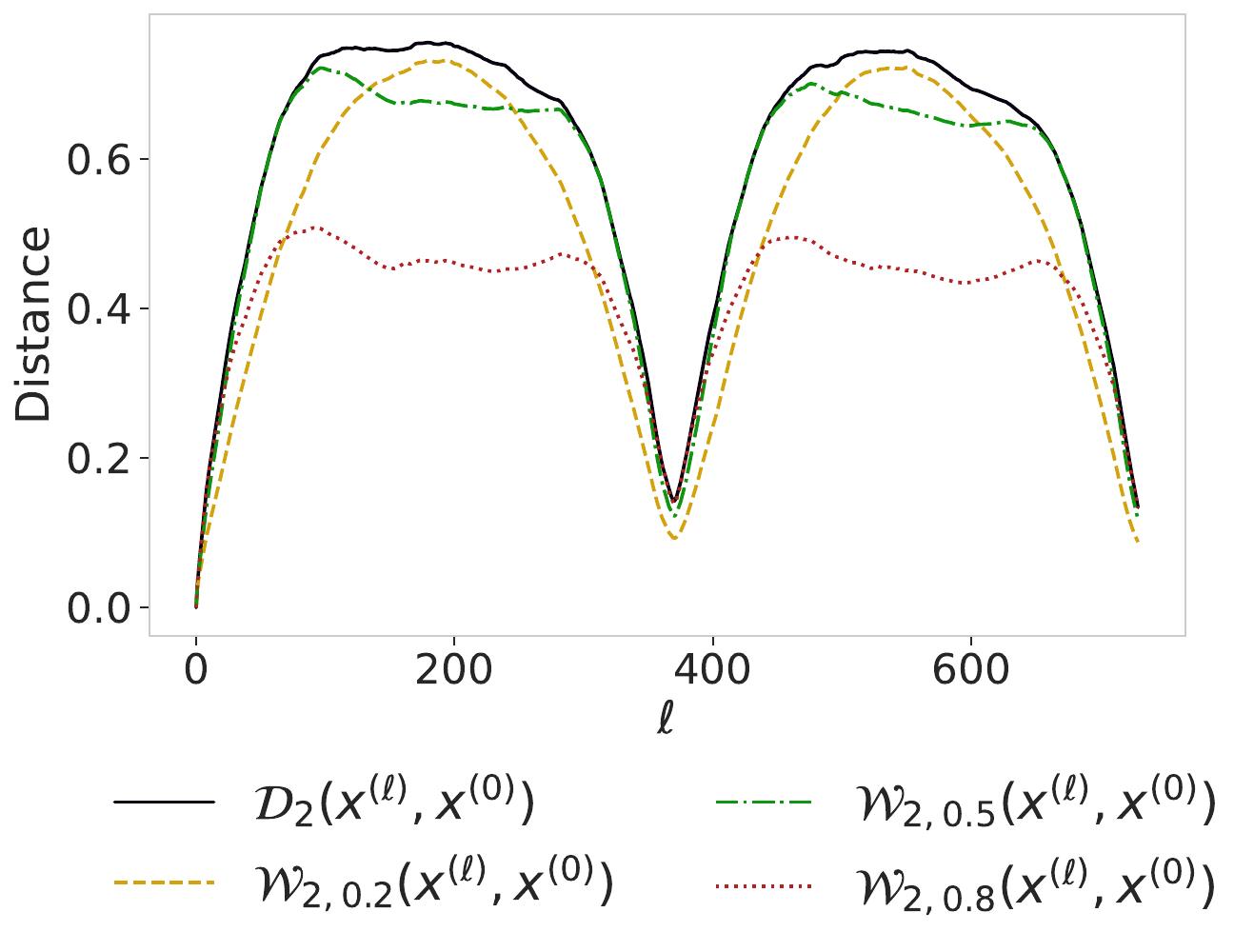}
        \label{fig:lag_distance_a}
    \end{subfigure}
    \begin{subfigure}[b]{0.495\textwidth}
        \includegraphics[width=\textwidth]{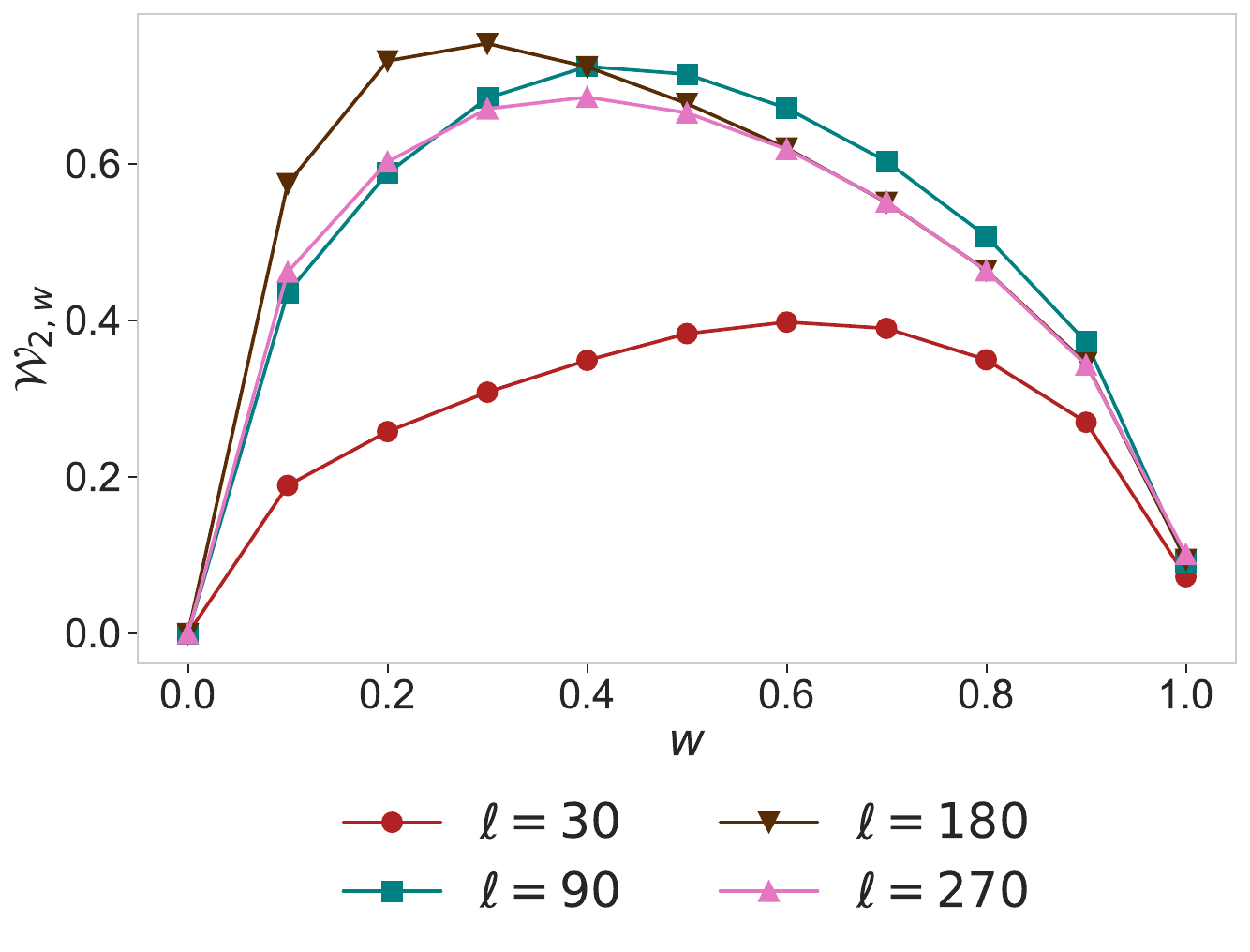}
        \label{fig:lag_distance_b}
    \end{subfigure}
    \vspace{-20pt}
    \caption{Left:  $\mathcal{D}_2$ and $\mathcal{W}_{2,w}$  between  \( x^{(0)}\) and \(x^{(\ell)}\). Right:  $\mathcal{W}_{2, w}$ metric against \(w\). }
    \label{fig:lag_distance}
\end{figure}

\paragraph{Analysis} In Figure \ref{fig:lag_distance} (left subplot), we report the values of  $\mathcal{D}_2$ and $\mathcal{W}_{2,w}$ between \(x^{(1)}\) and \(x^{(\ell)}\) against \(\ell\), where \(\ell=1,2,\dots,730\). All metrics capture the annual periodicity of the climate data, with dissimilarity values minimized near lags of 365 and 730, correctly reflecting seasonal similarity. However, for off-cycle lags such as when \(\ell\in[100,300]\) or \(\ell\in[450,650]\), the values of the metrics varies significantly. With a low weight on the temperature values (\(w=0.2\)), the $\mathcal{W}_{2, 0.2}$ produces an almost perfect sinusoid. Since the cost is defined as \(d_{2,w}^2 = w\|x-y\|^2_2 + (1-w)(t-s)^2 \), thus $\mathcal{W}_{2,0.2}$ recognizes the temporal shift but ignores the influence of the temperature variations. On the other hand, a higher weight \((w=0.8)\) returns a counter-intuitive result: \(\mathcal{W}_{2,0.8}(x^{(180)},x^{(0)})<\mathcal{W}_{2,0.8}(x^{(90)},x^{(0)}) \). This contradicts the natural expectation that a six-month shift should exceed a three-month one; similar behavior occurs for $w = 0.5$ and in the second year. In short, low $w$ ignores temperature information, while high $w$ distorts temporal relationships. In contrast, our proposed $\mathcal{D}_2$, interpretable as a maximum over all $\mathcal{W}_{2,w}$, achieves a robust balance between time and temperature, hence, retains temperature-specific fluctuations while clearly reflecting the seasonal cycle.

In Figure \ref{fig:lag_distance} (right subplot), we analyze the sensitivity of \(\mathcal{W}_{2,w}\) with respect to its weight parameter. We plot \(\mathcal{W}_{2,w}\) against \(w \in [0,1]\) for fixed shifts \(\ell \in \{30,90,180,270\}\) (days). When \(\ell\) is small (e.g., \(\ell=30\)), the $\mathcal{W}_{2,w}$ metric is relatively insensitive to \(w\), consistently yielding low dissimilarity values. For larger shifts (\(\ell=90,180,270\)), however, the choice of \(w\) becomes critical: \(\mathcal{W}_{2,w}\) varies substantially, and different values of \(w\) even change the ordering of which shift appears more significant. For example, when \(w\in [0.15,0.35]\), $\mathcal{W}_{2,w}$ considers \(\ell=180\) to be more significant, while for \(w\in [0.4,0.6]\), it considers \(\ell=90\) to be more significant. By contrast, our TiOT framework resolves this ambiguity by selecting the maximizer \(w\) for each pair of series, thereby providing a robust, parameter-free, and worst-case measure that is both consistent and reliable.

\subsection{Numerical performance of entropic TiOT}
\label{subsec: runtime_experiment}
In this section, we provide numerical experiments to validate our theoretical and computational findings. 

\begin{figure}[htb]
    \centering
        \begin{subfigure}[b]{0.495\textwidth}
        \includegraphics[width=\textwidth]{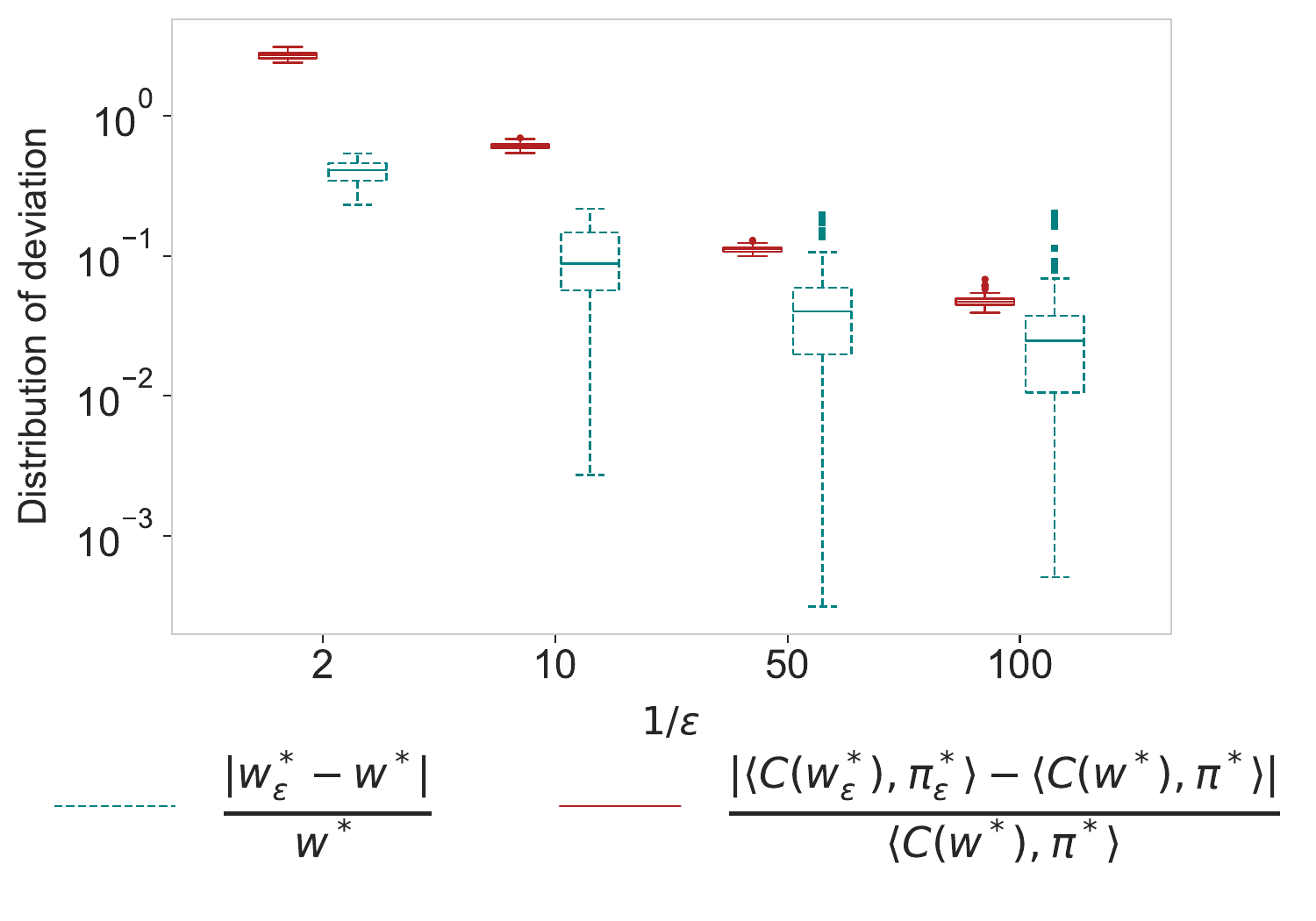}
    \end{subfigure}
        \begin{subfigure}[b]{0.495\textwidth}
        \includegraphics[width=\textwidth]{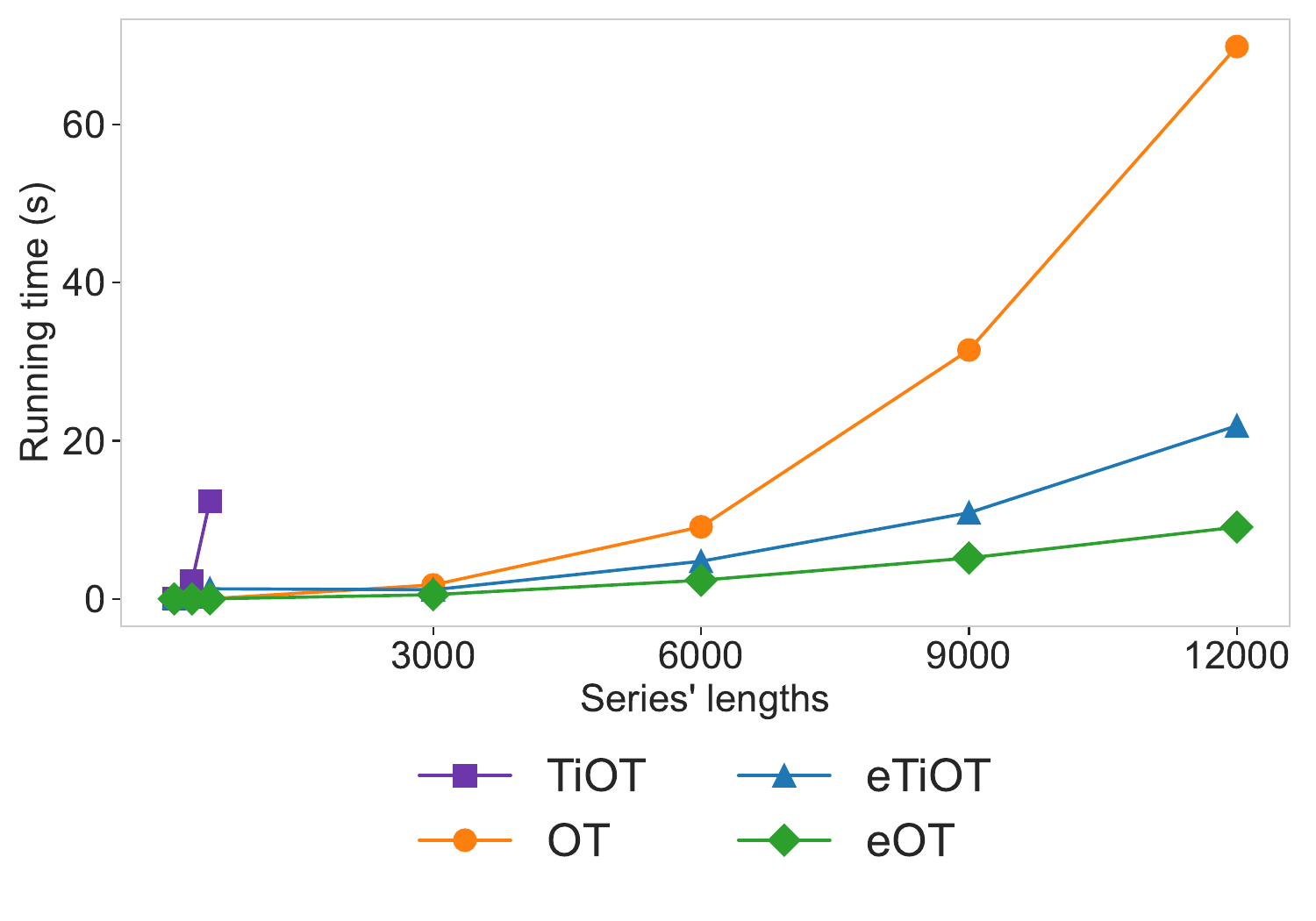}
    \end{subfigure}
    \vspace{-20pt}
 \caption{Left: Deviation of TiOT and eTiOT. Right: Computational comparison.}
 \label{fig:deviation_runtime}
\end{figure}

First, we demonstrate our theoretical result in Theorem~\ref{thm:convergence_eps}, that is, the convergence of the entropic TiOT (\namerefshort{def: eTiOT}{eTiOT}) to the exact \namerefshort{def: TiOT}{TiOT} as the regularization parameter \(\epsilon\) approaches zero.  For this experiment, we set  \(n= m = 100\) and generate two mixtures of Gaussians time series, similar to the procedure in Section~\ref{subsec: alignment_experiment}. We then solve the exact \namerefshort{def: TiOT}{TiOT} via its linear programming (LP) formulation (Appendix \ref{appen: LP_TiOT}) and its entropic regularized counterpart, \namerefshort{def: eTiOT}{eTiOT} (\ref{eq: eTiOT}) via Algorithm \ref{alg:HBCD} for a range of regularization parameters \(\epsilon = {1/2, 1/10, 1/50, 1/100} \). This process was repeated 100 times with different random seeds. Figure~\ref{fig:deviation_runtime} (left subplot) reports the distribution of deviations between the optimal objective values and the corresponding solutions of the two problems. One can observe that both differences decrease to 0 as \(\epsilon\rightarrow 0\). This result empirically verifies our theoretical analysis, showing that eTiOT serves as a reliable and accurate approximation of TiOT.

Second, we evaluate the computational complexity of our proposed algorithm against several relevant benchmarks. Specifically, we compare running time  of solving the eTiOT problem (\ref{eq: eTiOT}) via Algorithm~\ref{alg:HBCD} against its LP counterpart TiOT (\ref{eq: LP_TiOT}), as well as the standard OT (\ref{def: continuousTiOT}) and its entropic regularized version eOT, with the base metric $d_{2,w}$ and fixed parameter \(w=0.5\). For both eTiOT and eOT, we set the regularization parameter $\varepsilon = 0.1$ and the \(\operatorname{tolerance}=0.005\). Figure \ref{fig:deviation_runtime}
(right subplot) reports the running time required to solve each problem as the length of the time series \(n = m\) increases. As expected, the performance of TiOT, formulated as a large-scale LP, quickly becomes computationally prohibitive. The standard OT, solved by a standard Python Optimal Transport library (\cite{flamary2021pot, flamary2024pot}), also exhibits poor scaling. In contrast, Algorithm~\ref{alg:HBCD} shows a significant performance advantage, requiring only about 2–3 times the runtime of the highly efficient Sinkhorn algorithm for the classical eOT. 

\subsection{Experiments on 1NN classification}
\label{subsec: knn_experiment}
In this section, we assess the effectiveness of the distance defined by the TiOT problem for time series classification tasks. In particular, we compare the classification errors of the 1-nearest-neighbor algorithm using the following base metrics: Euclidean (ED), Dynamic Time Warping (DTW) with learned warping window \cite{Dau2019_UCR}, eTiOT (Algorithm \ref{alg:HBCD}), and eTAOT (Algorithm 1 in \cite{TAOT2020})\footnote{$\text{eTAOT}_{\omega}(\alpha, \beta) =  \langle C(\omega), \pi^*_\varepsilon \rangle$ where $ \pi^*_\varepsilon = \argmin_{\pi \in \Pi(\alpha, \beta)} \langle C(w), \pi \rangle + \mathbf{KL}(\pi | a \circ b)$ and $c_{ij}(\omega) = \|x-y\|^2_2 + \omega(t_i - s_j)^2$} on 15 datasets obtained from the UCR time series archive \cite{UCRArchive2018}. For ED and DTW, we adopt the classification errors provided by the benchmark website \cite{UCRArchive2018}.

First, the overall classification errors of ED, DTW, eTiOT, and eTAOT\((\omega)\) are presented in Table~\ref{table:error}. Both DTW and $\operatorname{eTAOT}(\omega)$ require a parameter controlling the strength of temporal constraints; these parameters are selected via leave-one-out cross-validation (LOOCV) in \cite{Dau2019_UCR, TAOT2020}, yielding $\text{learned\_w}$ for DTW and $w_{\text{grid}}$ for eTAOT. For the regularization parameter $\varepsilon$ of eTiOT and $\operatorname{eTAOT}(\omega)$, we perform 3-fold cross-validation over the grid 
$\{0.01,0.02, \dots, 0.1 \}$.

\begin{figure}[htb]
    \centering
    \includegraphics[width=0.9\linewidth]{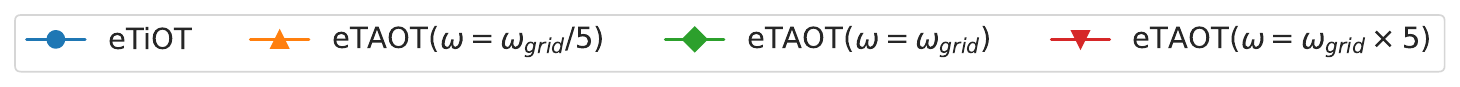}\\[3pt]
    \begin{subfigure}[b]{0.325\textwidth}
        \includegraphics[width=\textwidth]{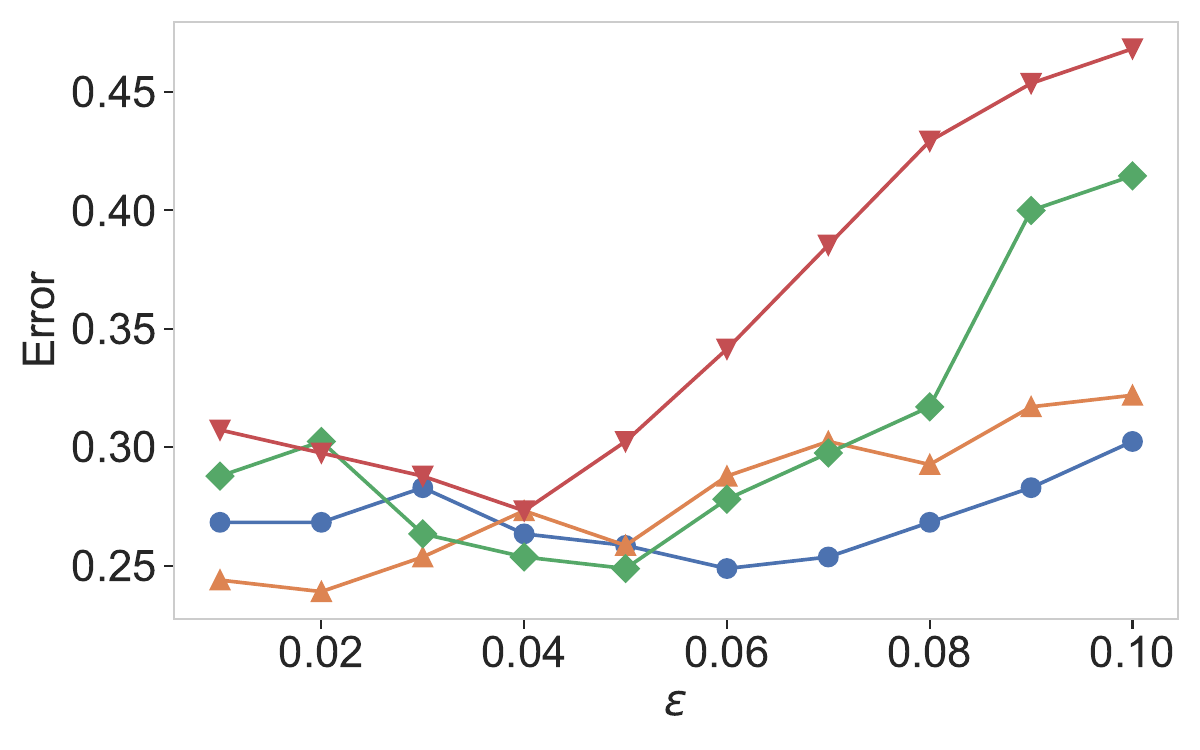}
        \caption{ProximalPhalanxTW}
    \end{subfigure}
    \begin{subfigure}[b]{0.325\textwidth}
        \includegraphics[width=\textwidth]{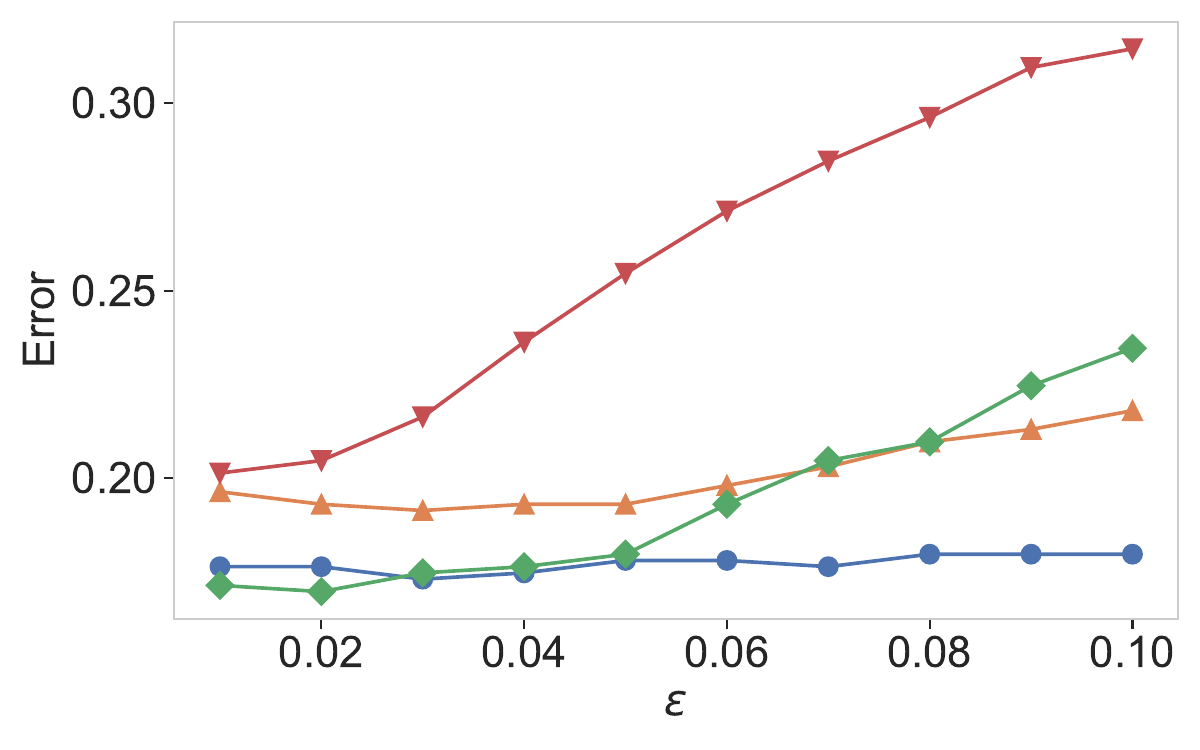}
        \caption{SonyAIBORobotSurface1}
    \end{subfigure}
    \begin{subfigure}[b]{0.325\textwidth}
        \includegraphics[width=\textwidth]{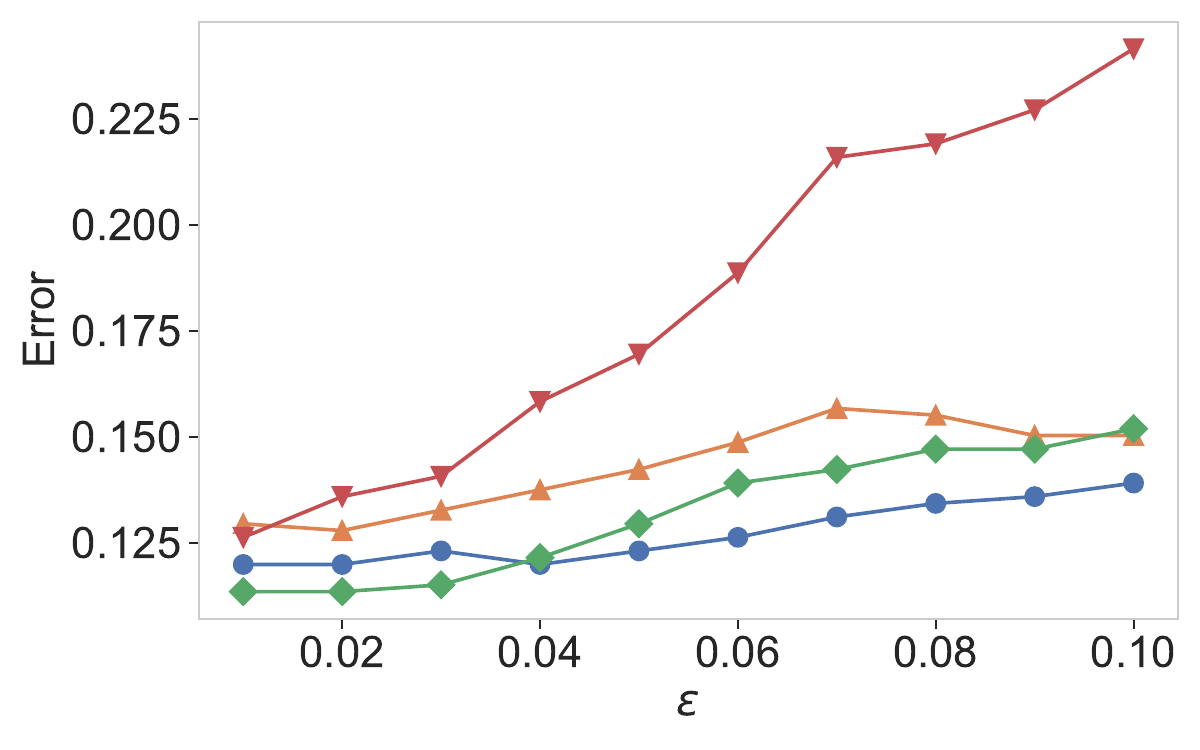}
        \caption{SwedishLeaf}
    \end{subfigure}
    \vspace{-10pt}
    \caption{Classification error of 1NN algorithm with $\varepsilon = 0.01, \dots, 0.1$.}
    \label{fig:errors}
\end{figure}

Second, we verify the robustness  of the eTiOT metric compared with $\operatorname{eTAOT}(\omega)$ for fixed $w$ by plotting classification errors across the  range $\varepsilon \in \{0.01, 0.02, \dots, 0.1 \}$. For $\operatorname{eTAOT}(\omega)$, we use the previously tuned $w_{\text{grid}}$ from LOOCV \cite{TAOT2020}; additionally, we include $\omega = \omega_{\text{grid}}/5$ and $\omega = 5\,\omega_{\text{grid}}$ to provide a more comprehensive comparison. Performance on 3 of the 15 datasets is shown 
in Figure~\ref{fig:errors}, while results for the remaining datasets are presented in Figure~\ref{fig: additional_errors} (Appendix \ref{appen: add_experiment}).

From Table~\ref{table:error}, the 1NN algorithm using the eTiOT metric achieves the lowest error on 10 of the 15 datasets. Specifically, its performance is better than or equal to ED on 12/15 datasets, DTW on 12/15 datasets, and eTAOT on 11/15 datasets. These results demonstrate both the effectiveness and the stability of eTiOT, highlighting its reliability as a dissimilarity measure for time series. Notably, eTiOT is the only metric in this experiment that does not require cross-validation to tune the strength of temporal constraints, relying instead on an adaptive mechanism. The robustness of eTiOT is further illustrated in Figure~\ref{fig:errors}, where the line representing eTiOT errors generally exhibits a slower increase as $\varepsilon$ grows from 0.01 to 0.1. It is well known that a larger regularization parameter typically reduces the running time of the entropic-regularized optimal transport problem. Therefore, this behavior not only provides additional evidence of eTiOT's robustness but also justifies the use of larger $\varepsilon$ values for faster computations.

\begin{table}[h]
\caption{Comparison of classification errors for \namerefshort{def: eTiOT}{eTiOT}, Euclidean, DTW, and eTAOT.}
\label{table:error}
\centering
\footnotesize
\begin{tabular}{lcc|cccc}
\hline
\textbf{dataset}      & \textbf{\makecell{train/test}} & \textbf{length} & \multicolumn{4}{c}{\textbf{error rates}}                          \\ \cline{4-7} 
                      & \textbf{size}                               &                 & ED             & DTW            & eTAOT          & eTiOT          \\ \hline
Adiac                 & 390/391                                     & 176             & 0.389          & 0.391          & 0.327          & \textbf{0.297} \\
ArrowHead             & 36/175                                      & 251             & \textbf{0.200} & \textbf{0.200} & 0.263          & 0.251          \\
BirdChicken           & 20/20                                       & 512             & 0.450          & 0.300          & 0.300          & \textbf{0.200} \\
CBF                   & 30/900                                      & 128             & 0.148          & \textbf{0.004} & 0.011          & \textbf{0.004} \\
DistalPhalanxOAG      & 400/139                                     & 80              & 0.374          & 0.374          & 0.324          & \textbf{0.317} \\
DistalPhalanxOC       & 600/276                                     & 80              & 0.283          & 0.275          & 0.261          & \textbf{0.257} \\
DistalPhalanxTW       & 600/276                                     & 80              & 0.367          & 0.367          & 0.367          & \textbf{0.353} \\
Ham                   & 109/105                                     & 431             & 0.400          & 0.400          & 0.390          & \textbf{0.362} \\
MiddlePhalanxOAG      & 400/154                                     & 80              & 0.481          & 0.481          & 0.487          & \textbf{0.468} \\
MiddlePhalanxOC       & 600/291                                     & 80              & \textbf{0.234} & \textbf{0.234} & 0.258          & 0.268          \\
MiddlePhalanxTW       & 600/291                                     & 80              & 0.487          & 0.494          & 0.455          & \textbf{0.422} \\
ProximalPhalanxOC     & 400/205                                     & 80              & \textbf{0.192} & 0.209          & 0.199          & 0.206          \\
ProximalPhalanxTW     & 400/205                                     & 80              & 0.293          & \textbf{0.244} & 0.249          & 0.268          \\
SonyAIBORobotSurface1 & 20/601                                      & 70              & 0.305          & 0.305          & 0.234          & \textbf{0.180} \\
SwedishLeaf           & 500/625                                     & 128             & 0.211          & 0.154          & \textbf{0.114} & 0.120          \\ \hline
\end{tabular}
\end{table}

\section{Conclusions}
\label{sec: conclusions}
In this work, we introduce \textit{Time-integrated Optimal Transport} (TiOT), a new framework for comparing time series. TiOT automatically balances temporal and feature information, thereby eliminating the need for manual parameter tuning compared to other measures. Moreover, we show that TiOT defines a proper metric that preserves fundamental properties of Wasserstein spaces.

We further develop an entropic regularized variant, eTiOT, and prove that it serves as a reliable approximation of TiOT. To solve eTiOT efficiently, 
we propose a Block Coordinate Descent (BCD) algorithm and provide a rigorous convergence analysis. Extensive experiments on synthetic and real-world datasets demonstrate the practical effectiveness and computational efficiency of our approach. Finally, TiOT offers a generalizable theoretical and practical foundation for defining a robust, adaptively tuned, weighted Euclidean distance between arbitrary data vectors.

\bibliographystyle{plain}
\bibliography{references}

\appendix

\section{Proof of the existence of optimal solution of \eqref{def: continuousTiOT}}\label{appen: proof_existence}
\begin{align*}
    \mathcal{D}_p(\alpha,\beta) 
    &= \max\limits_{w \in [0,1]} \left[ \min\limits_{\pi \in \Pi(\alpha, \beta)} 
    \int_{\mathbb{R}^{d+1} \times \mathbb{R}^{d+1}} d_{p,w}((x,t),(y,s))^p \, d\pi((x,t),(y,s)) \right]^{1/p} 
\end{align*}
By \cite[Theorem~1.7]{Samtambrogio2015}, the inner minimization admits an optimal solution. To prove the existence of an optimal $w^*$, we invoke the Weierstrass theorem. Since $[0,1]$ is compact, it remains to show that $\mathcal{T}_p(w) = \min_{\pi \in \Pi(\alpha, \beta)} \int d_{p,w}^p \, d\pi$ is upper semi-continuous with respect to $w \in [0,1]$.

To this end, we show the openness of the preimage $\mathcal{T}_p^{-1}((-\infty, a)) = \{w \in [0,1]: \mathcal{T}_p(w) < a \}$. Assume that $\mathcal{T}_p^{-1}((-\infty, a)) \neq \emptyset$, take an arbitrary $\bar{w} \in \mathcal{T}_p^{-1}((-\infty, a))$. Thus, there exists $\bar{\pi}$ such that $\int d_{p,\bar{w}}^p \, d\bar{\pi} < a$. Denote \( J_{\bar{\pi}} : [0,1] \to \mathbb{R} \) by \( J_{\bar{\pi}}(w) = \int d_{p,w}^p \, d\bar{\pi} \). Since $J_{\bar{\pi}}(w)$ is continuous with respect to $w$, the set $J_{\bar{\pi}}^{-1}((-\infty, a)) = \{w \in [0,1]: J_{\bar{\pi}}(w) < a \}$ is open. It follows that there exists $r > 0$ such that the ball $B(\bar{w}, r) \subset J_{\bar{\pi}}^{-1}((-\infty, a))$. Therefore, for all $w \in B(\bar{w}, r)$, we have $\mathcal{T}_p(w) \le J_{\bar{\pi}}(w) < a$, which implies $B(\bar{w}, r) \subset \mathcal{T}_p^{-1}((-\infty, a))$. Hence $\mathcal{T}_p^{-1}((-\infty, a))$ is open, and consequently, $\mathcal{T}_p(w)$ is upper semi-continuous.

\section{Proof of lemma \ref{lem: boundedness}}\label{appen: proof_bounded}
First, we have the bounds $C_{j\ell}(w^{k-1}) \geq C_{i\ell}(w^{k-1}) - 2\|C(w^{k-1})\|_\infty$  for all $i,j = 1,\dots,m,\, k = 1,\dots,n$, which gives
\begin{align*}
    u_i^k - u_j^{k} &=  - \varepsilon \log \big\langle e^{\frac{-C_{i \cdot}(w^{k-1})}{\varepsilon}}, e^{\frac{v^{k-1}}{\varepsilon}} \circ b \big\rangle + \varepsilon \log \big\langle e^{\frac{-C_{j \cdot}(w^{k-1})}{\varepsilon}}, e^{\frac{v^{k-1}}{\varepsilon}} \circ b \big\rangle \\
   &\leq -\varepsilon \log \big\langle  e^{\frac{-C_{i\cdot}(w^{k-1})}{\varepsilon}}, e^{\frac{v^{k-1}}{\varepsilon}} \circ b \big\rangle + \varepsilon \log \big\langle e^{\frac{-C_{i\cdot}(w^{k-1})}{\varepsilon} + \frac{2\|C(w^{k-1})\|_\infty}{\varepsilon}}, e^{\frac{v^{k-1}}{\varepsilon}} \circ b \big\rangle \\
   &= 2\|C(w^{k-1})\|_\infty
\end{align*}
where $C_{i\cdot}(w^{t-1})$ denotes the $i$-th row of $C(w^{t-1})$.
Combining this inequality and the normalization property, for any \(j=1,\dots,m\),  we have that 
\(- u_j^k = \sum_{i = 1}^m (u_i^k - u_j^k)a_i \leq 2\|C(w^{k-1})\|_\infty\)  
which implies 
\begin{equation*}
    u_j^k \geq - 2\|C(w^{k-1})\|_\infty.
\end{equation*}
Similarly, since \(u_i^k = \sum_{j = 1}^m (u_i^k - u_j^k)a_j \leq 2\|C(w^{k-1})\|_\infty\) for any \(i=1,\dots,m\), we have 
\begin{equation*}
    u_i^k \leq 2\|C(w^{k-1})\|_\infty.
\end{equation*}
Combining above two inequalities, we have $ \|u^k\|_\infty \leq 2\|C(w^{k-1})\|_\infty \leq 2\|C\|_\infty$. 
Now as $e^{\frac{-\|C(w^{k-1})\|_\infty}{\varepsilon}} \leq e^\frac{- C_{\ell i}(w^{k-1})}{\varepsilon} $ for any $\ell\in[m]$ and $i\in[n]$, we have
for all $i\in[n]$, 
\begin{align*}
    v_i^k \!= - \varepsilon \log \sum_{\ell = 1}^m 
    e^{\frac{-C_{\ell i}(w^{k-1})}{\varepsilon}} e^\frac{u_\ell^k}{\varepsilon} a_\ell 
    \leq -\varepsilon \log( e^\frac{-\|C(w^{k-1})\|_\infty - 2\|C(w^{k-1})\|_\infty}{\varepsilon} ) \!=\! 3 \|C(w^{k-1})\|_\infty.
\end{align*}
Applying an analogous argument and using the bounds $e^\frac{- C_{\ell i}(w^{k-1})}{\varepsilon}  \leq e^{\frac{\|C(w^{k-1})\|_\infty}{\varepsilon}}$, we have
$$v_i^k \geq -3\|C(w^{k-1})\|_\infty\quad\forall\; i\in[n].
$$
Therefore, we get $\|v^k_i\|_\infty \leq 3\|C(w^{k-1})\|_\infty \leq 3\|C\|_\infty$.

\section{Proof of lemma \ref{lem: block_lipschitz}}\label{appen: proof_block_lipschitz}
In the later part of this section, we shall frequently use the smoothness property of the exponential function over a bounded region; that is, given $M > 0$, for any $(a,b) \in [-M, M]\times [-M,M]$, we have
\begin{equation}
\vcenter{\hbox{$ |e^b - e^a| \leq e^M |b - a|.  $}}
\label{eq: exp: lipschitz}
\end{equation}
From this, we obtain
\begin{equation*}
   \begin{array}{lll}
    \!\! \text{LHS of \eqref{eq: block_lipschitz}} \!\! \!& = \! \!\!& \! \left|\sum_{i = 1}^m \sum_{j = 1}^n -\Tilde{C}_{ij}\exp(\frac{u_i^k + v_j^k}{\varepsilon} ) \left(\exp(\frac{-C_{ij}(w)}{\varepsilon}) - \exp(\frac{-C_{ij}(w')}{\varepsilon}) \right) a_i b_j \right|\\
    \! \!\!& \leq \! \!\!& \!(\|\Tilde{C}\|_\infty^2/\varepsilon)  \sum_{i = 1}^m \sum_{j = 1}^n \exp(\frac{u_i^k + v_j^k}{\varepsilon} ) \exp(\frac{\|C\|_\infty}{\varepsilon}) |w - w'| a_i b_j \\[5pt]
    \! \!\!& \leq \! \!\!& \!(\|\Tilde{C}\|_\infty^2/\varepsilon) \exp(6\|C\|_\infty / \varepsilon) |w - w'|,
\end{array} 
\end{equation*}
where $\Tilde{C}_{ij} = ||x_i - y_j ||^p_p - |t_i - s_j|^p$. The first inequality uses the inequality in \eqref{eq: exp: lipschitz}, and the definition of $\|\Tilde{C}\|_\infty = \max\{ |\Tilde{C}_{ij}| : i \in [m], j \in [n]\}$. The second inequality invokes the bounds from lemma \ref{lem: boundedness}.

\section{Proof of lemma \ref{lem: sufficient_decrease}}\label{appen: proof_sufficient_decrease}
In this proof, we shall use the strong convexity of the exponential function in a bounded region. Given $M > 0$, we have that for any $ (a,b) \in [-M, M] \times [-M, M]$,
\begin{align}
e^b - e^a - e^a(b - a) &\geq e^{-M}(b - a)^2/2, \label{eq: exp_strongconvex}
\end{align}
Let $\Delta_uF \!=\! F(u^k, v^k, w^k) \!-\! F(u^{k+1}, v^k, w^k)$, $\Delta_v F\!=\! F(u^{k+1}, v^k, w^k) \!-\! F(u^{k+1}, v^{k+1}, w^k)$ and $\Delta_w F =F(u^{k+1}, v^{k+1}, w^k) - F(u^{k+1}, v^{k+1}, w^{k+1})$. We have
\begin{equation*}
    F(u^k, v^k, w^k) - F(u^{k+1}, v^{k+1}, w^{k+1}) = \Delta_u F + \Delta_v F + \Delta_w F.
\end{equation*}
Next we will evaluate each of these terms separately. 
\begin{align}
    \Delta_u F &= \varepsilon \sum_{i= 1}^{m}\sum_{j=1}^n \Big(e^\frac{u_i^k}{\varepsilon} - e^\frac{u^{k+1}_i}{\varepsilon}\Big)e^\frac{v_j^k}{\varepsilon}e^\frac{-c_{ij}(w^k)}{\varepsilon}a_i b_j \label{eq: u1} \\
    &\geq \sum_{i=1}^{m}\sum_{j=1}^n (u_i^k - u_i^{k+1})e^\frac{u_i^{k+1}}{\varepsilon}e^\frac{v_j^k}{\varepsilon}e^\frac{-c_{ij}(w^k)}{\varepsilon}a_i b_j + \kappa \sum_{i = 1}^m (u_i^k - u_i^{k+1})^2  a_i \label{eq: u2}\\
    &= 
   e^\frac{\lambda^k}{\varepsilon}\sum_{i=1}^{m}(u_i^k - u_i^{k+1}) a_i
    + \kappa \sum_{i = 1}^m (u_i^k - u_i^{k+1})^2 a_i \label{eq: u3} \\
    &= \kappa  \|u^k - u^{k+1}\|^2_{L^2(\alpha)}, \label{eq: u4}
\end{align}
where \(\kappa \coloneqq \frac{e^\frac{-6\|C\|_\infty}{\varepsilon}}{2\varepsilon}\).
In the above, inequality~\eqref{eq: u2} is obtained by first
applying the local strong convexity~\eqref{eq: exp_strongconvex} with
$M = 2\|C\|_\infty/\varepsilon$ and then applying the bound
$\exp({\frac{v_j^k - c_{ij}(w^k)}{\varepsilon}}) \geq \exp({-\frac{4\|C\|_\infty}{\varepsilon}})$ to the second term. The equality (\ref{eq: u3}) is obtained by substituting the update expression of $u_i^{k+1}$ into $\exp({\frac{u_i^{k+1}}{\varepsilon}}$) in the first term of (\ref{eq: u2}).
The equality  (\ref{eq: u4}) comes from the fact that $(u^k - u^{k+1})^\top a = 0$ due to the normalization. 

Although $v^k$ has no normalization property, following analogous arguments, we still obtain
\begin{equation*}
    \Delta_v F \geq \kappa  \|v^k - v^{k+1}\|^2_{L^2(\beta)}.
\end{equation*}

In order to evaluate the last term, first we recall the following basic property of the projection operator onto a closed convex set $D \subset \mathbf{R}$:
\begin{equation}
\langle \mathbf{y} - P_D(\mathbf{y}), \mathbf{x} - P_D(\mathbf{y}) \rangle \leq 0 \quad \text{ for any } \mathbf{x} \in D, \mathbf{y} \in \mathbb{R}.
\end{equation}
Applying this property with $\mathbf{x} = w^k, \mathbf{y} = w^k - \eta \nabla_w F(u^{k+1}, v^{k+1}, w^k)$, we get
\begin{equation}
\left\langle \nabla_w F(u^{k+1}, v^{k+1}, w^k), w^k - w^{k+1} \right\rangle 
\geq \frac{1}{\eta} |w^k - w^{k+1}|^2.\label{eq: proj_property}
\end{equation}
Combining $\eta = (\varepsilon / \|\Tilde{C}\|^2_\infty) \exp({\frac{-6\|C\|_\infty}{\varepsilon}})$ with lemma~3.2, we obtain the Lipschitz continuity with constant $1/\eta$ of $\nabla_w F(u^k, v^k, w)$ over the compact set $[0,1]$ for all $k \geq 1$. Therefore, by invoking the  descent lemma \cite{Beck13_BCGD}, we have
\[
\Delta_w F \geq  - \left\langle \nabla_w F(u^{k+1}, v^{k+1}, w^k), w^{k+1} - w^k \right\rangle 
- \frac{1}{2\eta} |w^k - w^{k+1}|^2. 
\]
 Combining this with (\ref{eq: proj_property})  implies that
\[
\Delta_w F 
\geq \frac{1}{2\eta} |w^k - w^{k+1}|^2 = \tau |w^k -w^{k+1}|^2,
\]
where $\tau =  \|\Tilde{C}\|_\infty^2 \exp{ (\frac{ 6\|C\|_\infty}{\varepsilon})} / 2\varepsilon$.

\section{Proof of lemma \ref{lem: solution_boundedness}}\label{appen: proof_solution_boundedness}
    By Theorem \ref{thm: convergence_xi} and the equivalence of $\|\cdot\|_2$ and $\|\cdot\|_\infty$, we have $\inf_{\xi \in \mathcal{\Tilde{S}}}\|\xi^k - \xi\|_\infty \to 0$ as $k \to \infty$.
    By the definition of infimum, for each $k \geq 1$, we can choose $\xi^*_k \in \mathcal{\Tilde{S}}$ such that $\|\xi^k - \xi^*_k\|_\infty \leq 1/k + \inf_{\xi \in \mathcal{\Tilde{S}}}\|\xi^k - \xi\|_\infty$. Take the limit $k \to \infty$, we deduce $\|\xi^k - \xi^*_k\|_\infty \to 0$. Next, we consider
    \begin{equation}
    \label{eq: boundedsolution1}
        \begin{array}{lll}
            \|\xi^*_k\|_\infty  & \leq & \|\xi^k - \xi^*_k\|_\infty + \|\xi^k\|_\infty
             \; \leq \; \|\xi^k - \xi^*_k\|_\infty + \max\{1, 3\|C\|_\infty \},
        \end{array}
    \end{equation}
    where the first inequality uses the triangle inequality and the second inequality uses the bounds of $u^k, v^k, w^k$ from Lemma \ref{lem: boundedness}. Since $\{\xi^k - \xi^*_k\}$ converges, it is bounded. Therefore, from \eqref{eq: boundedsolution1}, we have that $\{\xi^*_k\}$ is bounded. Moreover, $\mathcal{\Tilde{S}}$ is closed, as it is the preimage of a closed set. Hence, there exists a subsequence $\{\xi^*_{k_i}\}$ which converges to $\xi^* = (u^*, v^*, w^*) \in \mathcal{\Tilde{S}}$. Using \eqref{eq: boundedsolution1}, we have
    \begin{equation*}
        \begin{array}{lll}
            \|u^*_{k_i}\|_\infty & \leq & \|u^{k_i} - u^*_{k_i}\|_\infty + \|u^{k_i}\|_\infty
             \;\leq \;\|u^{k_i} - u^*_{k_i}\|_\infty + 2\|C\|_\infty.
        \end{array}
    \end{equation*}
Taking $i \to \infty$ on both sides, we obtain $\|u^*\| \leq 2\|C\|_\infty$. Similarly, we can derive that $\|v^*\|_\infty \leq 3\|C\|_\infty$.

\section{Proof of theorem \ref{thm: sublinear}}\label{appen: proof_sublinear}
 The proof proceeds in four steps.\textbf{ Step 1} upper bounds the optimality gap $F(\xi^{k+1}) - F^*$ by the norm of the gradient difference. \textbf{Step 2} exploits the locally Lipschitz property of the exponential function to bound the norm of the gradient difference by the successive changes of the iterates. \textbf{Step 3} applies the sufficient descent property to bound the  successive changes of the iterates by the successive decreases in the function value. \textbf{Step 4} utilizes the inequality established from Step 3 and \cite[Theorem~3.1(1)]{Hong2017BCD} to immediately obtain the sublinear convergence rate of~\eqref{eq: HBCD}.

\textbf{Step 1:} let $\xi^* = (u^*, v^*, w^*)$ be an optimal solution of \eqref{eq: dual_eTiOT} satisfying $\|u^*\|_\infty \leq 2\|C\|_\infty$, $\|v^*\|_\infty \leq 3\|C\|_\infty$, its existence follows from Lemma \ref{lem: solution_boundedness}. Denote $\xi^{k+1/3} = (u^{k+1}, v^k, w^k)$ and $\xi^{k+2/3} = (u^{k+1}, v^{k+1}, w^k)$. By the convexity of \( F(\cdot) \), we have
\begin{equation}
\label{eq: sublinear_proof1}
    \begin{array}{lll}
          F(\xi^{k+1}) - F^* &\leq& \langle \nabla F(\xi^{k+1}), \xi^{k+1} - \xi^* \rangle \\
         &=& \langle \nabla_u F(\xi^{k+1}), u^{k+1} - u^* \rangle + \langle \nabla_v F(\xi^{k+1}), v^{k+1} - v^* \rangle \\[3pt]
         & & + \nabla_w F(\xi^{k+1})(w^{k+1} - w^*).
    \end{array}
\end{equation}
Next, we bound the terms on the RHS of \eqref{eq: sublinear_proof1} 
separately. First, we have that
\begin{equation*}
\begin{array}{@{\;}l@{\;}l@{\;}l@{\;}}
     \langle \nabla_u F(\xi^{k+1}), u^{k+1} - u^* \rangle & = &  \langle \nabla_u F(\xi^{k+1})  -  \nabla_u F(\xi^{k+1/3}) + \nabla_u F(\xi^{k+1/3}), u^{k+1} - u^* \rangle  \\[3pt]
    & \leq & \langle  \nabla_u F(\xi^{k+1})  -  \nabla_u F(\xi^{k+1/3}), u^{k+1} - u^* \rangle \\[3pt]
    & \leq & \| \ \nabla_u F(\xi^{k+1})  -  \nabla_u F(\xi^{k+1/3}) \|_2 \|u^{k+1} - u^* \|_2 \\[3pt]
    & \leq & \sqrt{m} \| \nabla_u F(\xi^{k+1})  -  \nabla_u F(\xi^{k+1/3}) \|_2 \|u^{k+1} - u^* \|_\infty \\[3pt]
    & \leq & 4\sqrt{m} \|C\|_\infty  \| \nabla_u F(\xi^{k+1})  -  \nabla_u F(\xi^{k+1/3}) \|_2.
\end{array} 
\end{equation*}
Here the first inequality is established by using the first-order optimality condition for $u^{k+1}$ with respect to $\varphi^k_u$, which gives $\langle \nabla_u F(\xi^{k+1/3}), u^{k+1} - u^* \rangle < 0$. The second inequality leverages the Cauchy–Schwarz inequality, the third follows from the equivalence of norms, and the last relies on the triangle inequality and the bounds of $\|u^k\|_\infty$ and $\|u^*\|_\infty$.

Similarly, utilizing the optimality of $v^{k+1}$, and $w^{k+1}$ of the upper bound functions $\varphi_v^k(v) = F(u^{k+1}, v, w^k)$ and $\varphi_w^k(w) = F(\xi^{k+2/3}) + \nabla_w F(\xi^{k+2/3})(w-w^k) + \frac{1}{2\eta}(w-w^k)^2$, we obtain
\begin{equation*}
    \begin{array}{lll}
         \langle \nabla_v F(\xi^{k+1}), v^{k+1} - v^* \rangle & \leq & 6\sqrt{n}\|C\|_\infty \|\nabla_v F(\xi^{k+1}) - \nabla_v F(\xi^{k + 2/3}) \|_2, \\[5pt]
         \nabla_w F(\xi^{k+1})(w^{k+1} - w^*) & \leq & | \nabla_w F(\xi^{k+1}) - \nabla \varphi_w^k (w^{k+1}) |.
    \end{array}
\end{equation*}

Denote $\mathbf{D}_{u}^k = \nabla_u F(\xi^{k+1})  -  \nabla_u F(\xi^{k+1/3})$, $\mathbf{D}_v^k = \nabla_v F(\xi^{k+1})  -  \nabla_v F(\xi^{k+2/3})$, and $\mathbf{D}_w^k = \nabla_w F(\xi^{k+1}) - \nabla \varphi_w^k (w^{k+1}) $, it follows from \eqref{eq: sublinear_proof1} and the above inequalities that
\begin{equation*}
    \begin{array}{lll}
        F(\xi^{k+1}) - F^* & \leq & 4\sqrt{m}\|C\|_\infty \|\mathbf{D}_u^k\|_2 + 6\sqrt{n}\|C\|_\infty\|\mathbf{D}_v^k\|_2 + |\mathbf{D}_w^k|, 
    \end{array}
\end{equation*}
which by Cauchy-Schwarz inequality implies
\begin{equation}
\label{eq: sublinear_proof2}
    \begin{array}{lll}
        (F(\xi^{k+1}) - F^*)^2 & \leq & 3 \left(16m\|C\|^2_\infty \|\mathbf{D}_u^k\|^2_2 + 36n\|C\|^2_\infty\|\mathbf{D}_v^k\|^2_2 + (\mathbf{D}_w^k)^2\right).
    \end{array}
\end{equation}

\textbf{Step 2: }we next proceed to bound $\|\mathbf{D}_{u}^k\|_2^2$ by component-wise analysis. Specifically, we examine each entry $\mathbf{D}_{u_i}^k$:
\begin{equation*}
    \begin{array}{lll}
         (\mathbf{D}_{u_i}^k)^2& \! = \! & \left[ \sum_{j = 1}^n \exp(\frac{u_i^{k+1}}{\varepsilon})\left( \exp(\frac{v_j^{k+1} - C_{ij}(w^{k+1})}{\varepsilon}) - \exp(\frac{v_j^{k} - C_{ij}(w^{k})}{\varepsilon})  \right)a_i b_j \right]^2 \\[3pt]
         & \! \leq \! & a_i^2  \sum_{j = 1}^n \exp(\frac{2u_i^{k+1}}{\varepsilon})\left( \exp(\frac{v_j^{k+1} - C_{ij}(w^{k+1})}{\varepsilon}) - \exp(\frac{v_j^{k} - C_{ij}(w^{k})}{\varepsilon})  \right)^2 b_j \\[3pt]
         &\! \leq \! & \frac{a_i^2}{\varepsilon^2} \exp(\frac{12\|C\|_\infty}{\varepsilon}) \left( \sum_{j = 1}^n (v_j^{k+1} - v_j^k + \Tilde{C}_{ij}(w^{k+1} - w^k))^2  b_j \right) \\[5pt]
         &\! \leq \!& \frac{2a_i^2}{\varepsilon^2}  \exp(\frac{12\|C\|_\infty}{\varepsilon}) \left( \sum_{j = 1}^n (v_j^{k+1} - v_j^k)^2 b_j + \Tilde{C}_{ij}^2 (w^{k+1} - w^k)^2 b_j \right) \\[5pt]
         &\! \leq \!& \frac{2a_i^2}{\varepsilon^2}  \exp(\frac{12\|C\|_\infty}{\varepsilon}) \left( \|v^{k+1} - v^k\|^2_{L^2(b)} + \|\Tilde{C}\|^2_\infty  (w^{k+1} - w^k)^2 \right),
    \end{array}
\end{equation*}
where $\Tilde{C}_{ij} = ||x_i - y_j ||^p_p - |t_i - s_j|^p$. The first inequality uses the Cauchy-Schwarz inequality and the property $\left(\sum_{j = 1}^n b_j \right) = 1$, the second leverages \eqref{eq: exp: lipschitz} and the bounds from lemma \ref{lem: boundedness}, the third again exploits the Cauchy-Schwarz inequality. Finally, the last inequality uses the definitions of $\|\cdot\|_{L(b)}$ and $\|\Tilde{C}\|_\infty$. Given that $\left(\sum_{i = 1}^m a_i^2 \right) \leq 1$, we deduce that
\begin{equation}
\begin{array}{lll}
\label{eq: D_u}
    \| \mathbf{D}^k_u \|_2^2  & \leq & \frac{2}{\varepsilon^2}   \exp(\frac{12\|C\|_\infty}{\varepsilon}) \left( \|v^{k+1} - v^k\|^2_{L^2(b)} + \|\Tilde{C}\|^2_\infty  (w^{k+1} - w^k)^2 \right).
\end{array}
\end{equation}
Using similar arguments, we obtain that
\begin{equation}
\label{eq: D_v}
    \begin{array}{lll}
        \|\mathbf{D}_v^k \|_2^2 & \leq & \frac{\|\Tilde{C}\|^2_\infty}{\varepsilon^2}\exp(\frac{12\|C\|_\infty}{\varepsilon}) (w^{k+1} - w^k)^2.
    \end{array}
\end{equation}
To complete Step 2, we evaluate
\begin{equation}
\label{eq: D_w}
    \begin{array}{lll}
        |\mathbf{D}^k_w| & = & |\nabla F_w(\xi^{k+1}) - \nabla F_w(\xi^{k+2/3}) - (1/\eta)(w^{k+1} - w^k)|  \\[3pt]
        & \leq & | \nabla F_w(\xi^{k+1}) - \nabla F_w(\xi^{k+2/3})| + \frac{\|\Tilde{C}\|_\infty^2}{\varepsilon}\exp(\frac{6\|C\|_\infty}{\varepsilon})|w^{k+1} - w^k| \\[3pt]
        & \leq & \sum_{i = 1}^m \sum_{j = 1}^n \|\Tilde{C}\|_\infty \exp(\frac{5\|C\|_\infty}{\varepsilon})\left|\exp(\frac{-C_{ij}(w^{k+1})}{\varepsilon}) - \exp(\frac{-C_{ij}(w^{k})}{\varepsilon})\right|a_i b_j \\[3pt]
        & &  + (\|\Tilde{C}\|_\infty^2/\varepsilon)\exp(\frac{6\|C\|_\infty}{\varepsilon})|w^{k+1} - w^k| \\[3pt]
        & \leq & (2\|\Tilde{C}\|_\infty^2/\varepsilon) \exp(\frac{6\|C\|_\infty}{\varepsilon})|w^{k+1} - w^k| \\[3pt]
        & \leq & (2\|\Tilde{C}\|_\infty \|C\|_\infty/\varepsilon) \exp(\frac{6\|C\|_\infty}{\varepsilon})|w^{k+1} - w^k|,
    \end{array}
\end{equation}
where the first inequality follows from the triangle inequality and the formula of the stepsize~$\eta$, the second uses the triangle inequality of the absolute value and the bounds of $u^{k+1}, v^{k+1}, \Tilde{C}_{ij}$, the third leverages \eqref{eq: exp: lipschitz}, and the last uses the fact that $\|\Tilde{C}\|_\infty \leq \|C(1)\|_\infty \leq \|C\|_\infty$.

Now, substitute \eqref{eq: D_u}, \eqref{eq: D_v}, \eqref{eq: D_w} into \eqref{eq: sublinear_proof2}, we arrive at
\begin{equation*}
\begin{array}{@{\:}r@{\:}c@{\:}l@{}}
 (F(\xi^{k+1}) - F^*)^2
 & \le &
 96m\rho \!\left(
      \|v^{k+1} - v^k\|_{L^2(b)}^2
      + \|\tilde C\|_\infty^2 (w^{k+1} - w^k)^2
 \right)
 \\[3pt]
 &&
 + \rho\!\left(
      108n\|\tilde C\|_\infty^2 (w^{k+1} - w^k)^2
      + 12\|\tilde C\|_\infty^2 (w^{k+1} - w^k)^2
 \right)
 \\[3pt]
 & \le &
\rho (96m + 108n + 12) \!\left(
      \|v^{k+1} - v^k\|_{L^2(b)}^2
      + \|\tilde C\|_\infty^2 (w^{k+1} - w^k)^2
 \right),
\end{array}
\end{equation*}
where $\rho = \frac{\|C\|_\infty^2 }{\varepsilon^2}\exp(\frac{12\|C\|_\infty}{\varepsilon})$.

\textbf{Step 3:} By using the sufficient descent property (Lemma \ref{lem: sufficient_decrease}) to upper bound the RHS of the above inequality, we obtain 
\begin{equation}
\label{eq: sublinear_proof3}
    \begin{array}{lll}
         (F(\xi^{k+1}) - F^*)^2 &\hspace{-0.8em} \leq \rho_1 \left(F(\xi^k) - F(\xi^{k+1})\right),
    \end{array}
\end{equation}
where $\rho_1 = (192m + 216n + 24)\frac{\|C\|^2_\infty}{\varepsilon}\exp(\frac{18\|C\|_\infty}{\varepsilon})$.

\textbf{Step 4:} From the analysis presented in \cite[Theorem 3.1(1)]{Hong2017BCD} and the inequality \eqref{eq: sublinear_proof3}, the proof is  completed.

\section{TiOT as a Linear Program}
\label{appen: LP_TiOT}
Following similar steps as in Section~\ref{sec: BCD}, we can rewrite the \namerefshort{def: TiOT}{TiOT} problem as
\begin{equation}
\min\limits \left\{ q^\top z \,\mid\,
Hz \leq r, w \in [0,1], \; z=[u;v;w] \in \mathbb{R}^{\,m+n+1}\right\} 
\label{eq: LP_TiOT}
\end{equation}
where 
$q = [a; b;0] \in \mathbb{R}^{\,m+n+1},$ the constraint matrix is
\[
H \in \mathbb{R}^{\,mn \times (m+n+1)}, \quad 
H_{ i + (j-1)n,:} = \bigl( -e_i^\top , \; -e_j^\top , \; (t_i - s_j)^2 - \|x_i - y_j\|^2 \bigr),
\]
and \(r_{i + (j-1)n} = (t_i - s_j)^2\), for \(i\in[m]\) and \(j\in[n]\).

\newpage
\section{Additional experiments on time series classification}
\label{appen: add_experiment}
\begin{figure}[h]
    \centering
    \begin{subfigure}[b]{0.3122\textwidth}
        \includegraphics[width=\textwidth]{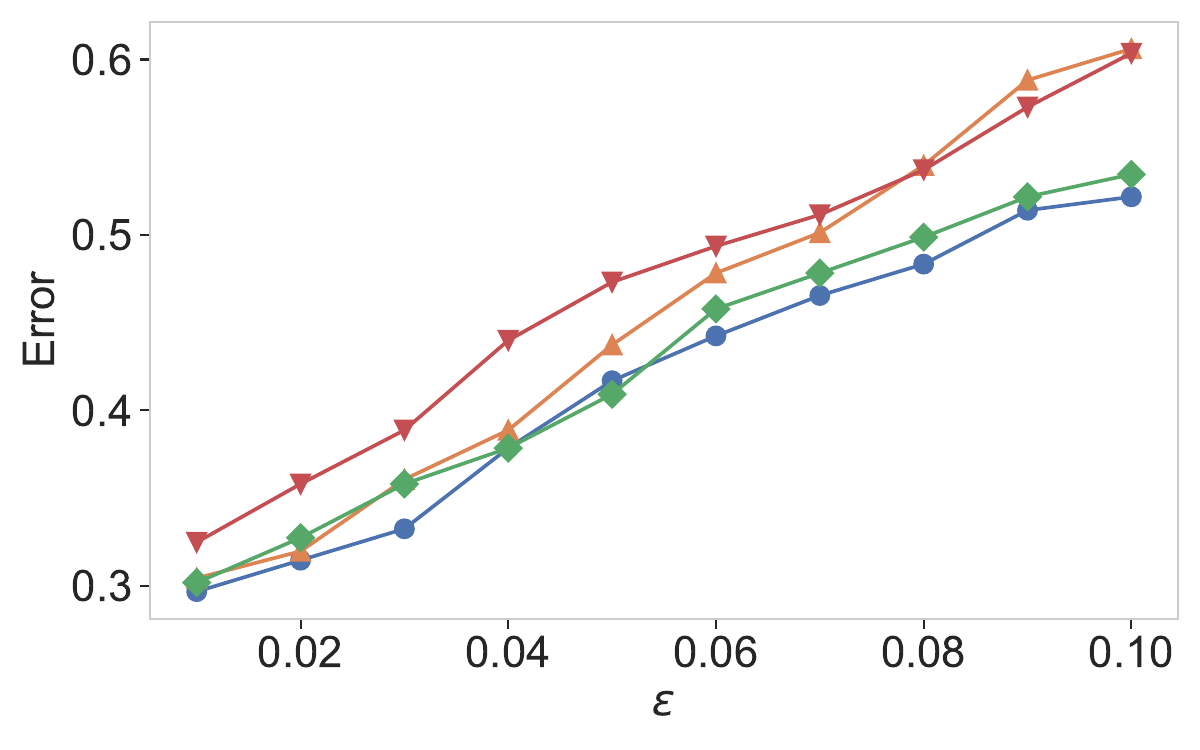}
        \caption{Adiac}
    \end{subfigure}
    \begin{subfigure}[b]{0.3122\textwidth}
        \includegraphics[width=\textwidth]{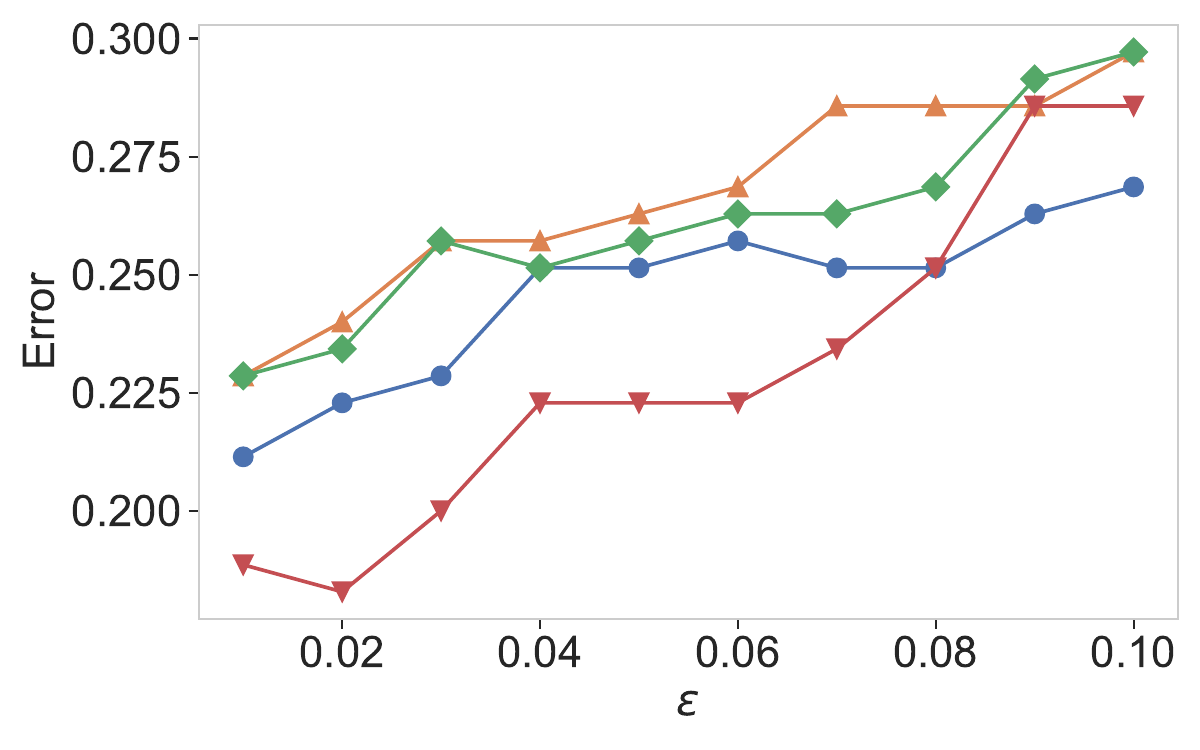}
        \caption{ArrowHead}
    \end{subfigure}
    \begin{subfigure}[b]{0.3122\textwidth}
        \includegraphics[width=\textwidth]{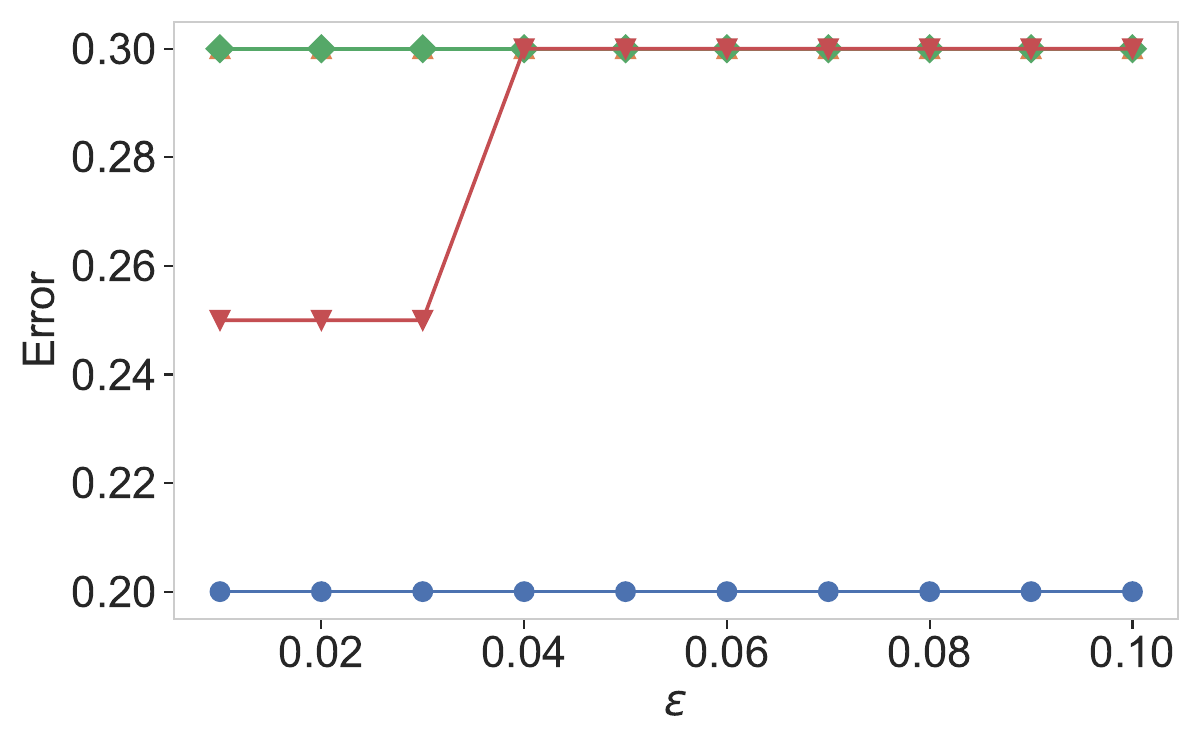}
        \caption{BirdChicken\protect\footnotemark}
        \label{fig:figure2}
    \end{subfigure} 
    \begin{subfigure}[b]{0.3122\textwidth}
        \includegraphics[width=\textwidth]{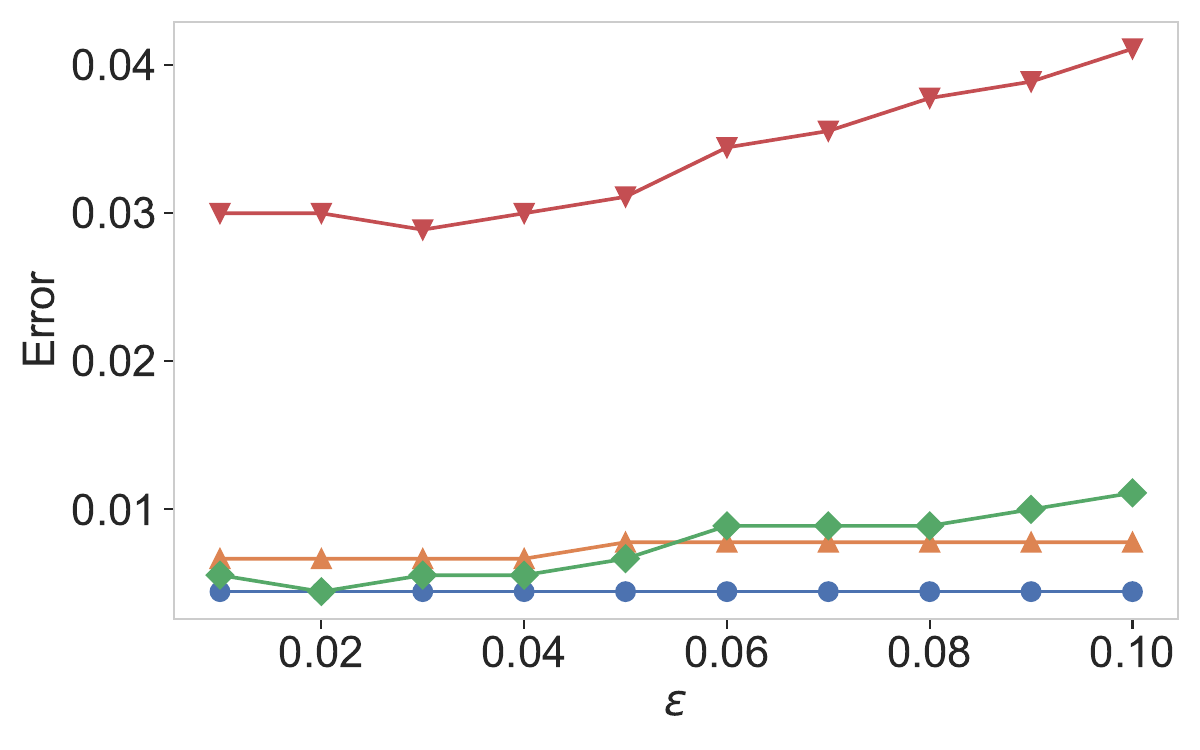}
        \caption{CBF}
        \label{fig:figure2}
    \end{subfigure}
        \begin{subfigure}[b]{0.3122\textwidth}
        \includegraphics[width=\textwidth]{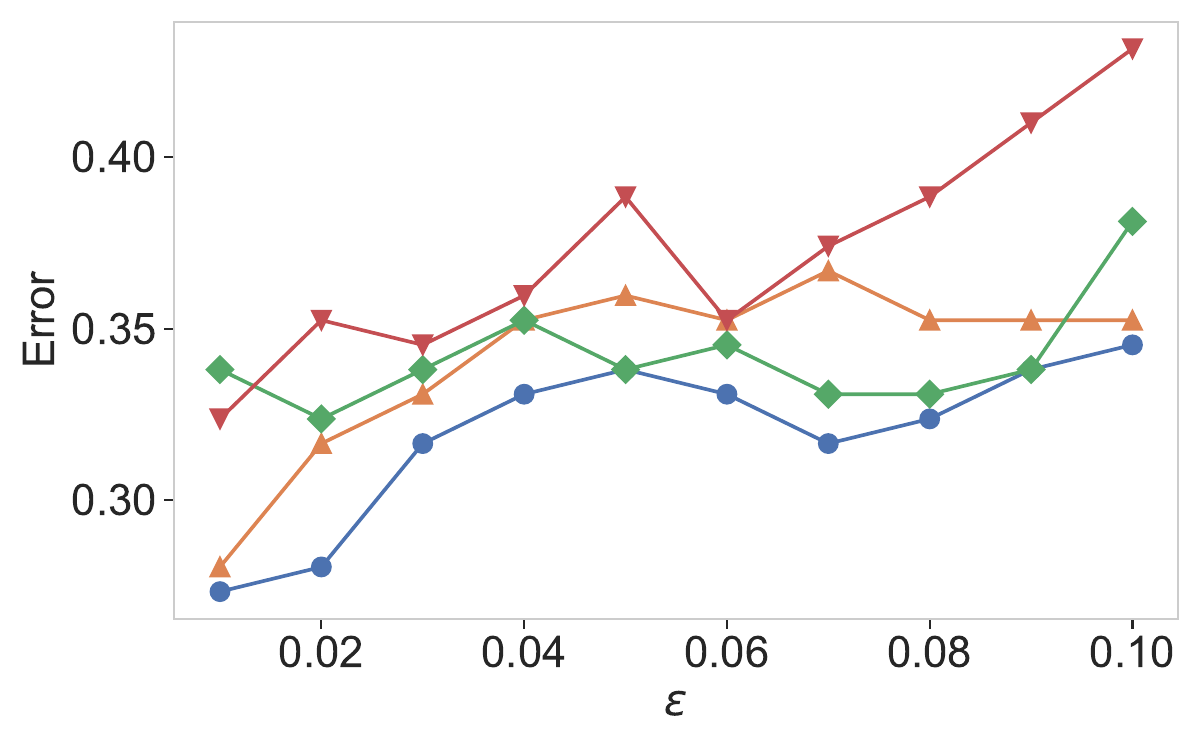}
        \caption{DistalPhalanxOAG}
    \end{subfigure}
        \begin{subfigure}[b]{0.3122\textwidth}
        \includegraphics[width=\textwidth]{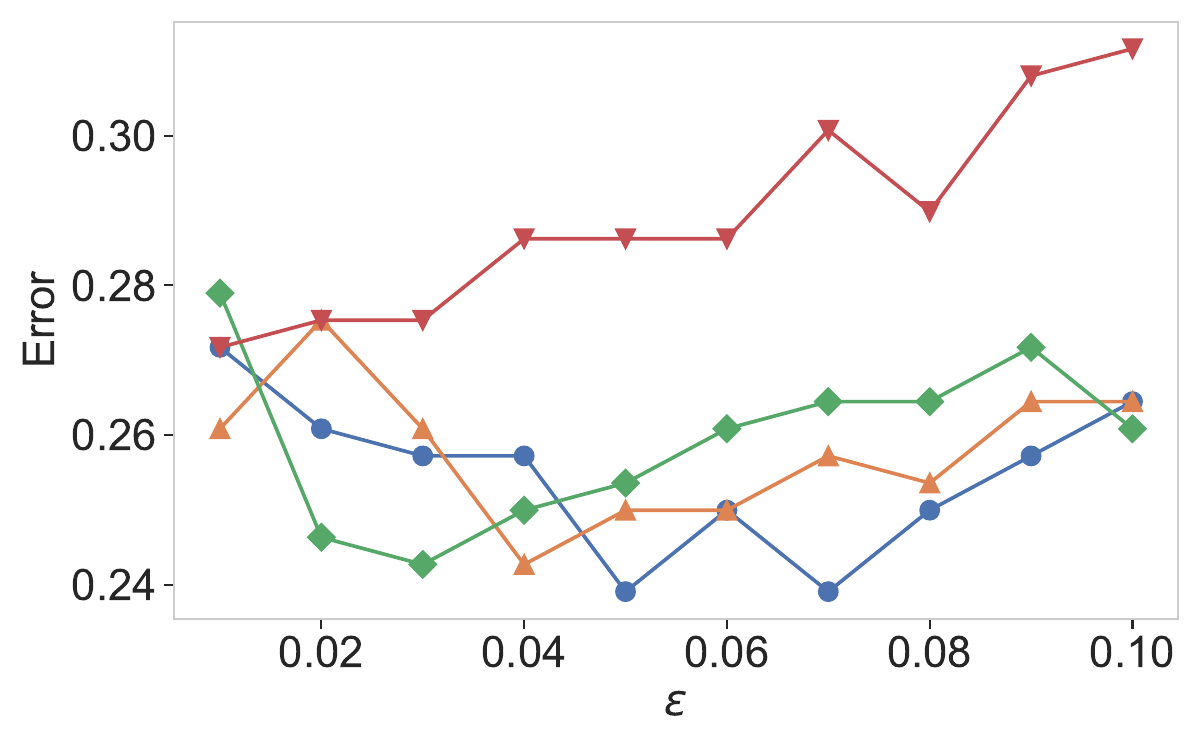}
        \caption{DistalPhalanxOC}
    \end{subfigure}
    \begin{subfigure}[b]{0.3122\textwidth}
        \includegraphics[width=\textwidth]{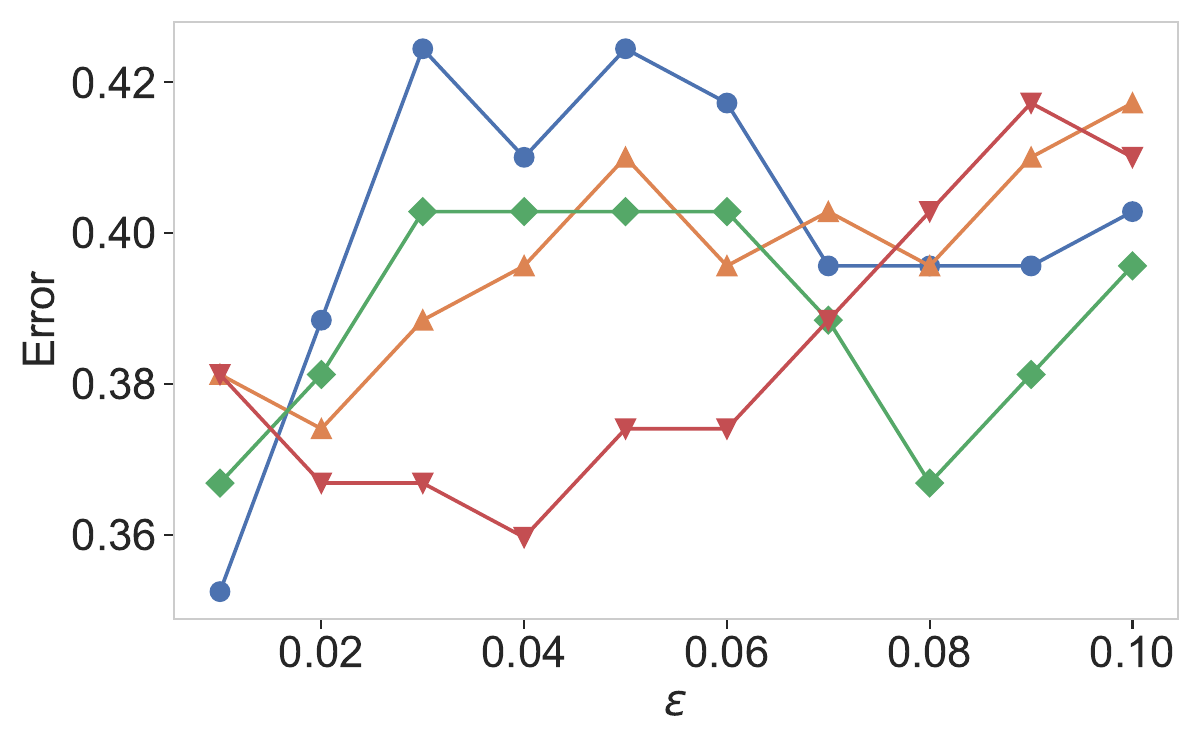}
        \caption{DistalPhalanxTW}
    \end{subfigure}
    \begin{subfigure}[b]{0.3122\textwidth}
        \includegraphics[width=\textwidth]{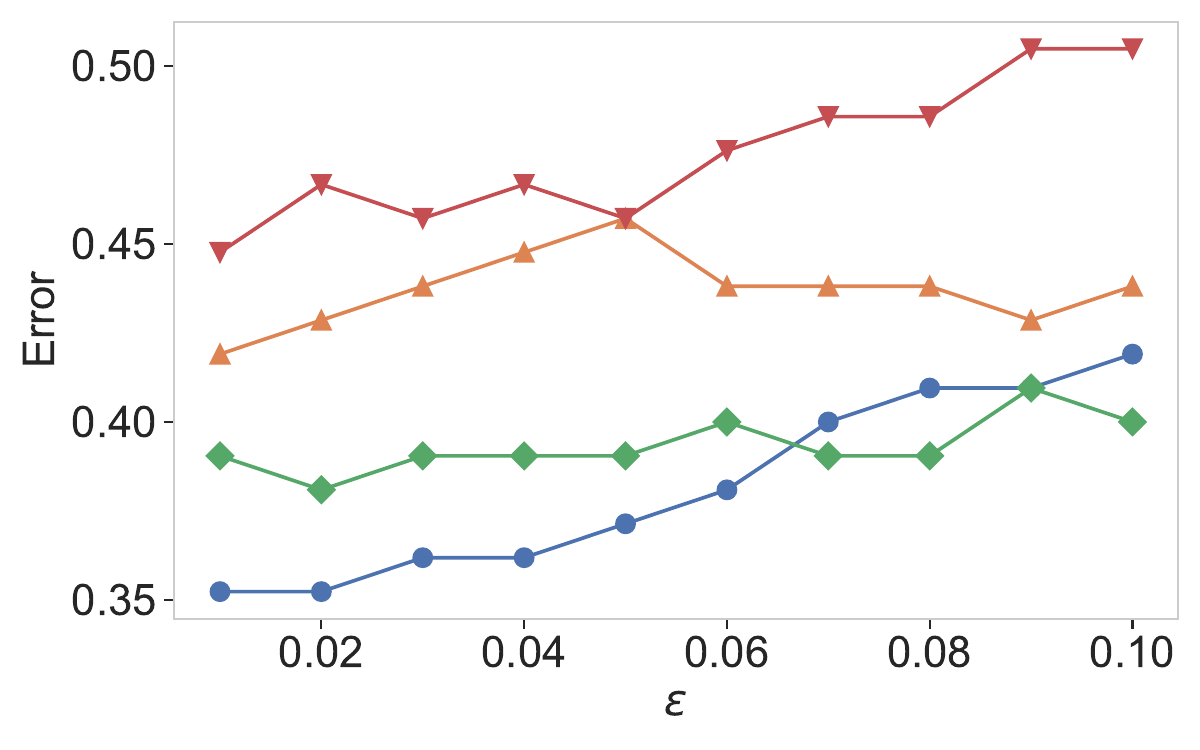}
        \caption{Ham}
    \end{subfigure}
    \begin{subfigure}[b]{0.3122\textwidth}
        \includegraphics[width=\textwidth]{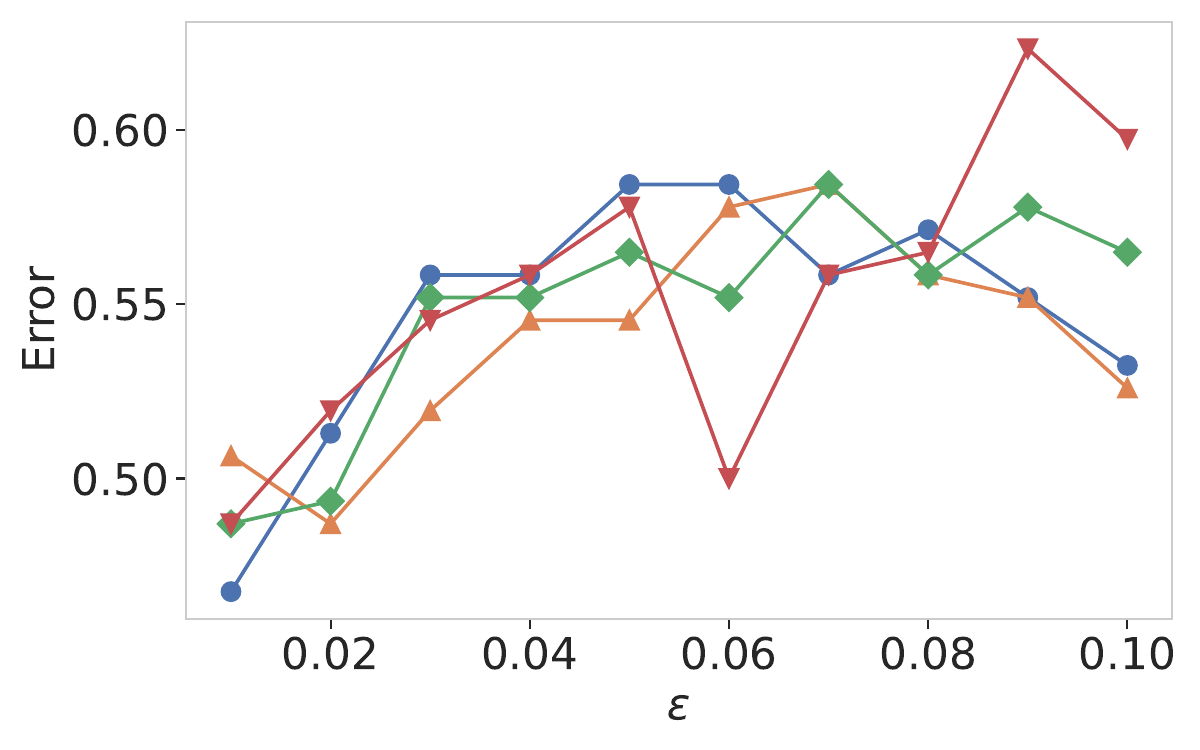}
        \caption{MiddlePhalanxOAG}
    \end{subfigure}
    \begin{subfigure}[b]{0.3122\textwidth}
        \includegraphics[width=\textwidth]{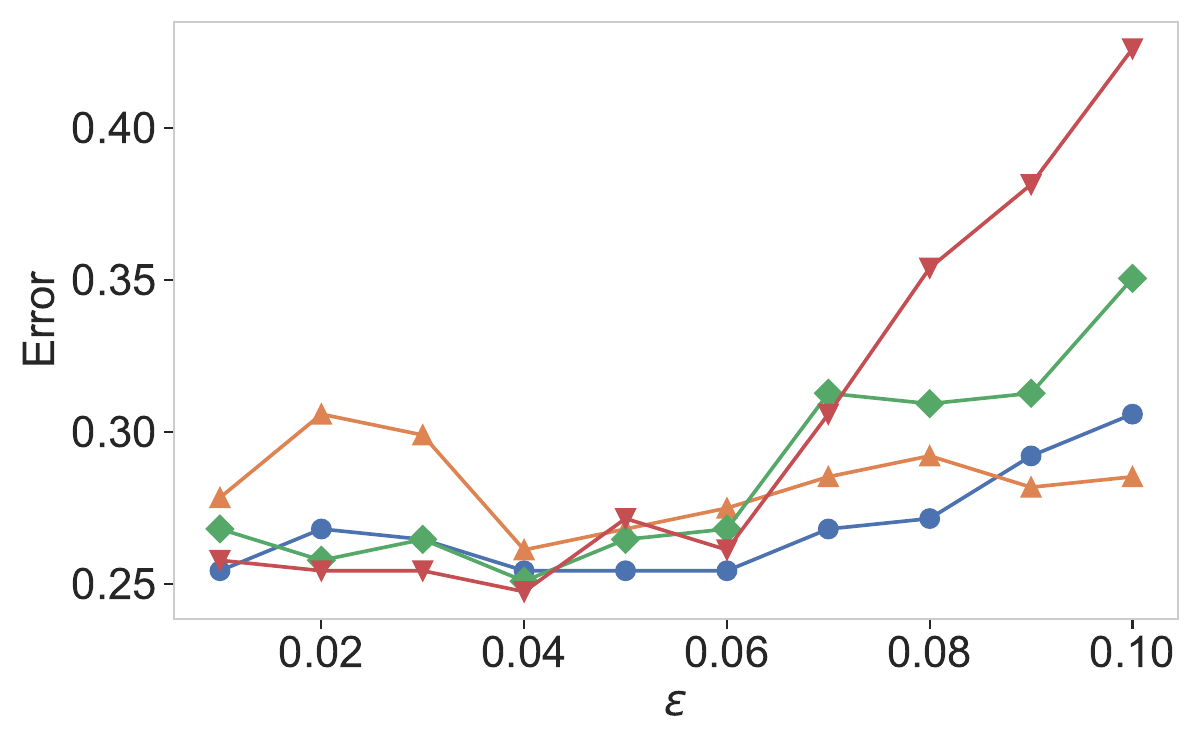}
        \caption{MiddlePhalanxOC}
    \end{subfigure}
        \begin{subfigure}[b]{0.3122\textwidth}
        \includegraphics[width=\textwidth]{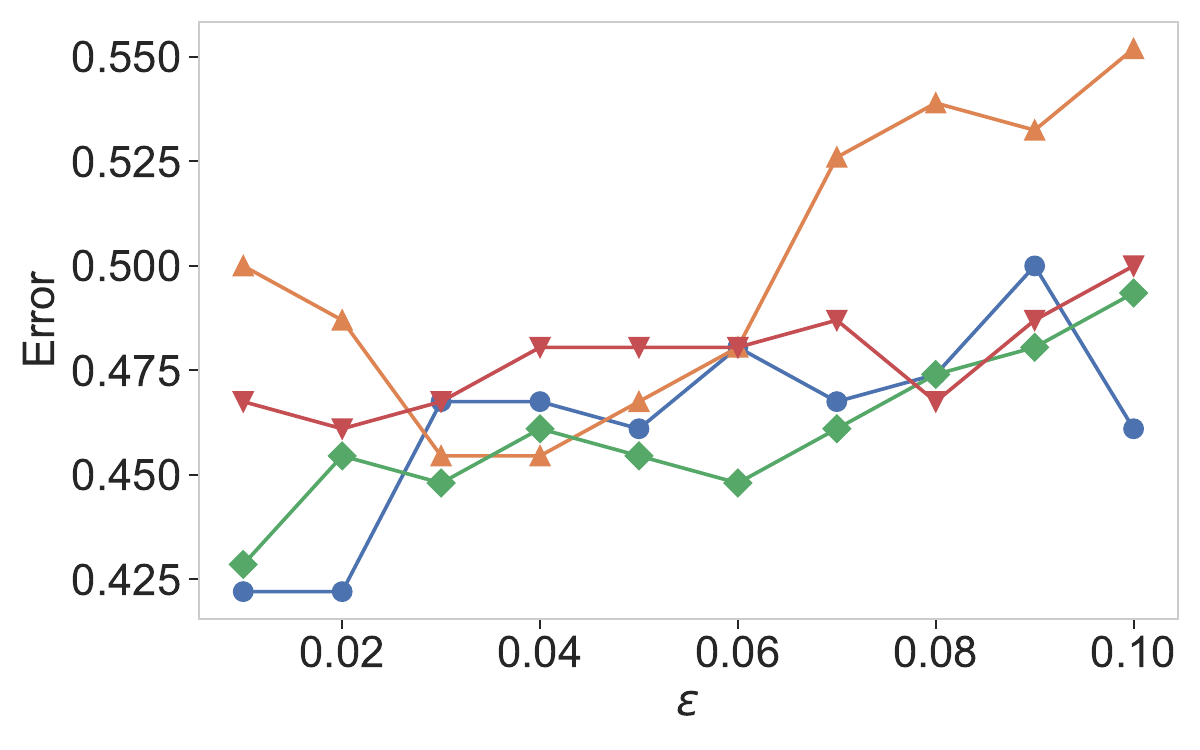}
        \caption{MiddlePhalanxTW}
    \end{subfigure}
        \begin{subfigure}[b]{0.3122\textwidth}
        \includegraphics[width=\textwidth]{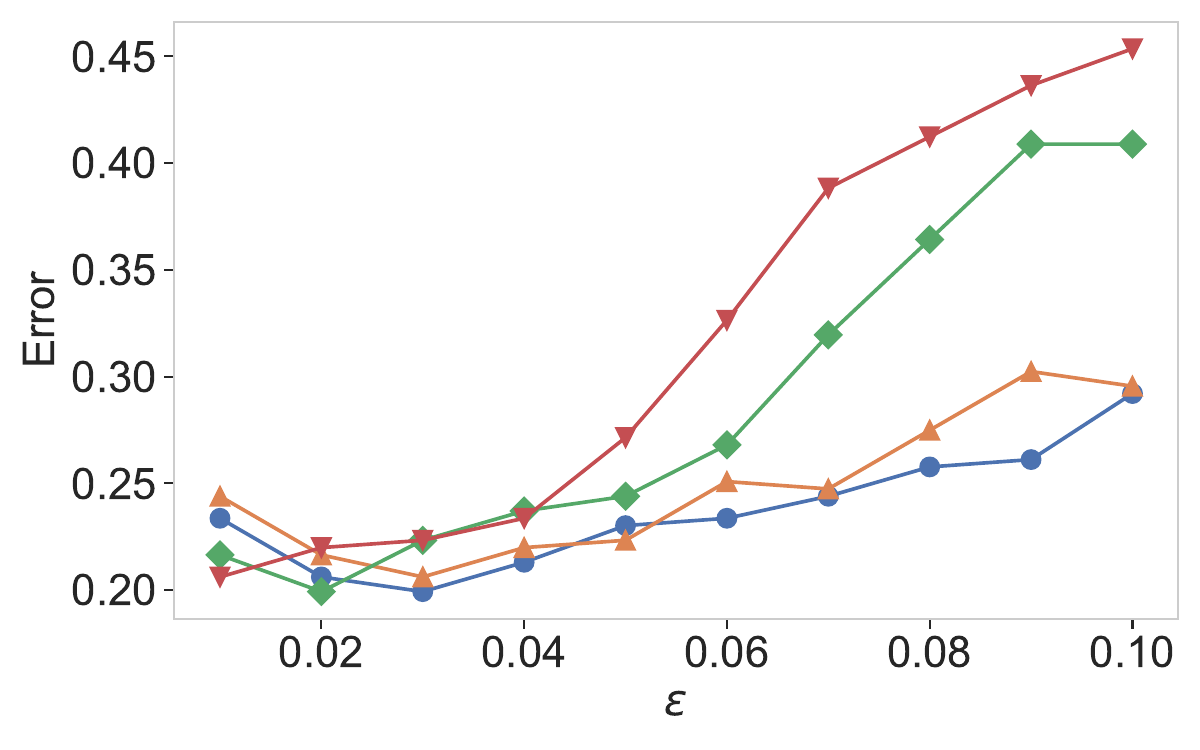}
        \caption{ProximalPhalanxOC}
    \end{subfigure}
    \includegraphics[width=0.9\linewidth]{Figures/legend_horizontal.pdf}\\[3pt]
    \caption{Classification error of 1NN algorithm with $\varepsilon = 0.01, \dots, 0.1$.}
    \label{fig: additional_errors}
\end{figure}
\footnotetext{For the BirdChicken dataset, due to the small test set,
$\operatorname{eTAOT}(w = w_{\text{grid}}/5)$ and
$\operatorname{eTAOT}(w = w_{\text{grid}})$ produce identical errors (0.30).
}

\end{document}